\documentclass[11pt,a4paper]{article}

\usepackage{amsmath,amssymb}
\newcommand{\e}{\varepsilon}
\newcommand{\D}{\Delta}

\newcommand{\va}{\varphi}
\newcommand{\n}{\nabla}
\newcommand{\N}{\frac{N}{2}}

\newcommand{\T}{\mathbb{T}}
\newcommand{\p}{\partial}

\newcommand{\R}{\mathbb{R}}
\newcommand{\NN}{\frac{N}{p}}

\newcommand{\de}{\delta}

\newtheorem{definition}{Definition}
\newtheorem{theorem}{Theorem}

\newtheorem{proposition}{Proposition}
\newtheorem{corollaire}{Corollary}

\newtheorem{remarka}{Remark}

\newtheorem{lemme}{Lemma}

\title{Existence of global strong solution for the compressible Navier-Stokes system and the Korteweg system in two-dimension}
\author{Boris Haspot \thanks{Ceremade UMR CNRS 7534
Universit\'e de Paris  Dauphine,
Place du MarŽchal DeLattre De Tassigny
75775 PARIS CEDEX 16 , haspot@ceremade.dauphine.fr 
}}
\date{}
\begin{document}
%\tableofcontents
\maketitle
\begin{abstract}
 This paper is dedicated to the study of viscous compressible barotropic fluids in dimension $N=2$. We address
the question of the global existence of strong solutions with large initial data for compressible Navier-Stokes system and Korteweg system. In the first case we are interested by slightly extending  a famous result due to V. A. Vaigant and A. V. Kazhikhov in \cite{VG} concerning the existence of global strong solution in dimension two for a suitable choice of viscosity coefficient ($\mu(\rho)=\mu>0$ and $\lambda(\rho)=\lambda\rho^{\beta}$ with $\beta>3$) in the torus. We are going to weaken the condition on $\beta$ by assuming only $\beta>2$ essentially by taking profit of commutator estimates introduced by Coifman et al in \cite{4M} and using a notion of \textit{effective velocity} as in \cite{VG}. In the second case we study the existence of global strong solution with large initial data in the sense of the scaling of the equations for Korteweg system with degenerate viscosity coefficient and with friction term. %We shall use the notion of quasi solution introduced in \cite{arxiv,cras1,cras2}.
% , by the fact that we choose initial data more general in $\widetilde{B}^{\NN-1,\NN}_{p,1}\times B^{\frac{N}{p_{1}}}_{p_{1},1}$ with $1\leq p, p_{1}<+\infty$ good chosen. Moreover we extend the class of viscosity coefficients.
\end{abstract}
\section{Introduction}The motion of a general barotropic compressible fluid with capillary tensor is described by the following system, which can be derived from a Cahn-Hilliard  free energy (see the
pioneering work by J.- E. Dunn and J. Serrin in \cite{3DS} and also in
\cite{3A,3C,3GP}):
\begin{equation}
\begin{aligned}
\begin{cases}
&\p_{t}\rho+{\rm div}(\rho u)=0,\\
&\p_{t}(\rho u)+{\rm div}(\rho u\otimes u)-{\rm div}(2\mu(\rho)D u)-\n\big(\lambda(\rho){\rm div}u\big)+\n P(\rho)+a\rho u={\rm div}K,
\end{cases}
\end{aligned}
\label{1}
\end{equation}
where ${\rm div}K$ is the capillary tensor which reads as follows:
\begin{equation}
{\rm div}K
=\n\big(\rho\kappa(\rho)\D\rho+\frac{1}{2}(\kappa(\rho)+\rho\kappa^{'}(\rho))|\n\rho|^{2}\big)
-{\rm div}\big(\kappa(\rho)\n\rho\otimes\n\rho\big).
\label{divK}
\end{equation}
The term
${\rm div}K$  allows to describe the variation of density at the interfaces between two phases, generally a mixture liquid-vapor. $P$ is a general increasing pressure term that we assume in the sequel under the form $P(\rho)=b\rho^{\gamma}$ with $b>0$ and $\gamma\geq1$, $a\rho u$ is a friction term with $a>0$ (see \cite{P}).
%In the sequel to simplify, $\kappa(\rho)$ will take the form $\kappa(\rho)=\rho^{\alpha}$
%with $\alpha\in\R$.\\
$D (u)=\frac{1}{2}[\n u+^{t}\n u]$ being the stress tensor, $\mu$ and $\lambda$ are the two Lam\'e viscosity coefficients depending on the density $\rho$ and satisfying:
$$\mu>0\;\;\mbox{and}\;\;2\mu+N\lambda\geq0.$$
In the present paper, we are interested in dealing with two different situations:
\begin{itemize}
\item The case of the compressible Navier-Stokes system where we assume no capillarity, $\kappa(\rho)=0$ and where $\mu(\rho)=1$ is a constant and $\lambda(\rho)=\lambda\rho^{\beta}$ with $\beta\geq 2$.
\item The case of Korteweg system with the viscosity coefficients and the capillarity coefficient such that:
$$\mu(\rho)=\mu\rho\;\;\lambda(\rho)=0\;\;\mbox{and}\;\;\kappa(\rho)=\frac{\kappa}{\rho}\;\;\mbox{with}\;\;\kappa>0,\;\mu>0.$$
\end{itemize}
In the first case we would like to extend the famous result of global strong solution in two dimension in the torus discovered by V. A. Vaigant and A. V. Kazhikhov in \cite{VG}. Indeed in \cite{VG}, the authors assume that $\lambda(\rho)=\lambda\rho^{\beta}$ with $\beta>3$ and $\lambda>0$ and $\mu(\rho)=1$. Let us emphasize that a such choice on the viscosity coefficients allows to exhibits two different phenomena; the first one concerns the notion of \textit{effective velocity} introduced by Lions in \cite{PL2}  which is crucial in this context for getting a priori estimates on the divergence and the rotational of the velocity; the second reason to choose such coefficient concerns the possibility to obtain $L^{\infty}_{T}(L^{p}(\mathbb{T}^{2}))$ estimates for any $p>1$ on the density and any $T>0$. Indeed the viscosity coefficient $\lambda(\rho)$ offers enough weight in order to derivate such estimates on the density (see the p1115 "Second a priori estimate for the density" in \cite{VG}). In the first part of this paper we wish to improve this result by assuming only $\lambda(\rho)=\lambda\rho^{\beta}$ with $\beta>2$ and $\lambda>0$. The key point will consist in using commutator estimates for dealing with term of the form $[R_{ij}, u_{j}](\rho u_{i})$ (we refer to  \cite{4M} for such estimates but also to \ref{L2} where Lions prove global existence of weak solution for compressible Navier-Stokes equations by introducing this kind of ingredient).\\
%In the second case we are interested in proving global existence of strong solution for the Korteweg system with friction term when:
%\begin{equation}
%\mu(\rho)=\mu\rho\;\;\lambda(\rho)=0\;\;\mbox{and}\;\;\kappa(\rho)=\frac{\kappa}{\rho}\;\;\mbox{with}\;\;\kappa>0,\;\mu>0;
%\label{cas2}
%\end{equation}
%and:
%\begin{equation}
%P(\rho)=\alpha\rho^{\gamma}\;\;\mbox{and}\;\;a=-\mu\alpha.
%\label{cas2a}
%\end{equation}
In the second case we are interested in proving global existence of strong solution for the Korteweg system with friction term when the physical coefficients verify:
\begin{equation}
\kappa(\rho)=\frac{\kappa}{h},\;\kappa=\mu^{2}\;\;\mbox{and}\;\;b=a\mu,
\label{hypothese}
\end{equation}
with $\mu>0$. This system without friction has been widely studied this last year in particular concerning the existence of global weak solution and global strong solution with small initial data. We refer in particular to the following works \cite{3BDL,3DD,Hprepa,Hprepa1,fH1,fH2,H1,J}.
Let us start with explaining a other notion of \textit{effective velocity} used in particular in \cite{Hprepa} which allows us to simplify the system (\ref{1}). Indeed by computation (see \cite{Hprepa}), we obtain the simplified system:
%\begin{equation}
%\begin{cases}
%\begin{aligned}
%&\p_{t}\rho+{\rem div}(\rho v)-\frac{\kappa}{\mu}\D\rho=0,\\
%&\rho\p_{t}v +\rho u\cdot\n v-\rm div(\mu\rho\,\n v)+\n P(\rho)=0,
%\end{aligned}
%\end{cases}
%\label{3systeme1}
%\end{equation}
%When $\alpha=0$ and $\kappa=\mu^{2}$, we obtain the following simplified model:
\begin{equation}
\begin{cases}
\begin{aligned}
&\p_{t}\rho+{\rm div}(\rho v)-\mu\D \rho=0,\\
&\rho\p_{t}v +\rho u\cdot\n v-\rm div(\mu \rho\,\n v)+a \,\rho\,v=0,
\end{aligned}
\end{cases}
\label{3systeme2}
\end{equation}
with $v=u+\mu\n\ln \rho$ the \textit{effective velocity}. For more details on the computation, we refer to \cite{Hprepa1}. When we write the system (\ref{3systeme2}) in function of the momentum $m=\rho v$, the system reads as follows:
\begin{equation}
\begin{cases}
\begin{aligned}
&\p_{t}\rho+{\rm div}m-\mu\D \rho=0,\\
&\p_{t}m +\rm div(\frac{m}{\rho}\otimes m)-\mu\D m+r\,m=0,
\end{aligned}
\end{cases}
\label{3system3}
\end{equation}
In particular we observe that $(\rho,-\mu\n\ln\rho)$ when:
$$\p_{t}\rho-\mu\D \rho=0,$$
is a particular global solution of the system (\ref{1}).
\begin{remarka}
Let us mention that we can choose initial density which admits  vacuum that we wish. In general it is not always possible to obtain global strong solution with initial density close from the vacuum.
\label{global}
\end{remarka}
%In the sequel we shall obtain two type of results, a first one where we assume that the friction coefficient verifies the assumption (\ref{hypothese}). 
In the sequel we will be interested in working around this global particular solution (see remark \ref{global}) in order to prove the existence of global strong solution with large initial data on the irrotational part.%\\
%In a second time, we would like to avoid this condition of friction by assuming that $a=0$. In this case we have no explicit solution as it is mentioned in the remark \ref{global}. In order to overcome this difficulty, we are going to introduce the notion of \textit{quasi solution} developed in \cite{arxiv,cras1,cras2} and to work around this last one.
%\begin{definition}
%We say that $(\rho,u)$ is a quasi solution of the system (\ref{1}) if $(\rho,u)$ verifies the following system almost everywhere:
%\begin{equation}
%\begin{aligned}
%\begin{cases}
%&\p_{t}\rho+{\rm div}(\rho u)=0,\\
%&\p_{t}(\rho u)+{\rm div}(\rho u\otimes u)-{\rm div}(\mu(\rho)\D u)-\n\big(\lambda(\rho){\rm div}u\big)={\rm div}K,
%\end{cases}
%\end{aligned}
%\label{2}
%\end{equation}
%\begin{equation}
%\end{equation}
%\end{definition}
%In this case we would like to mention the existence of global strong solution with large initial data for the scaling of the equations. More precisely it means that our smallness assumption can be chosen supercritical, it should allow us to exhibit large initial data for the energy space in dimension $N=2$ generating global strong solution. Let us recall that with this choice of viscosity coefficient the existence of global strong solution with large initial in two dimension remains open in general, in this sense this result give a first kind of answer to the problem. 
More precisely we shall wish obtain global strong solution in critical space for the scaling of the equations. Let us briefly recall the notion of invariance by scaling of the equation and by what we mean by supercritical smallness on the data. By critical, we mean that we want to solve the system (\ref{1}) in functional spaces with invariant norm by the natural changes of scales which leave (\ref{1}) invariant. More precisely in 
our  case,  the following transformation:
\begin{equation}
(\rho(t,x),u(t,x))\longrightarrow (\rho(l^{2}t,lx),lu(l^{2}t,lx)),\;\;\;l\in\R,
\label{1}
\end{equation}
verify this property, provided that the pressure term has been changed accordingly. In particular we can observe that $\dot{H}^{\N}\times \dot{H}^{\N-1}$ is a space invariant for the scaling of the equation, more generally such Besov spaces:
$$(\rho_{0}-\bar{\rho})\in B^{\NN}_{p,1},\;u_{0}\in B^{\frac{N}{p_{1}}-1}_{p_{1},1},$$
with $(p,p_{1})\in [1,+\infty[$ are also available. 
%\begin{definition}
%We say that the smallness assumption on the initial data is supercritical if $(\rho_{0}-1,u_{0})$ can be arbitrarily large in critical Besov space.
%\end{definition}
%\begin{remarka}
%Let us emphasize that generally all the result of global strong solution for the system (\ref{1}) involved a smallness condition in  critical Besov space.
%\end{remarka}
\subsection{Results}
Let us state the two main result of this paper. The first one is an improvement of the results of Vaigant and Kazhikhov \cite{VG}.
\begin{theorem}
\label{theo1}
Let us assume the following hypothesis on the viscosity coefficients:
$$\mu(\rho)=1\;\;\mbox{and}\;\;\lambda(\rho)=\lambda\rho^{\beta}\;\;\mbox{with}\;\beta>2.$$
Let $u_{0}\in H^{2}(\T^{2})$, $\rho_{0}\in W^{1,q}(\T^{2})$ with $q>2$ and:
$$0<c\leq \rho_{0}(x)\leq m<+\infty\;\;\forall x\in\T^{2},$$
then it exists a unique global strong solution to (\ref{1}) such that:
$$u\in W^{2,1}_{2}(Q_{T})\;\;\mbox{and}\;\;\rho\in W^{1,1}_{q,\infty}(Q_{T})\;\;\forall T>0,$$
with $Q_{T}=(0,T)\times\mathbb{T}^{2}$.
\end{theorem}
\begin{remarka}
As mentioned above, the main point compared with the result of \cite{VG} concerns the fact that we can improve the range of $\beta$ by assuming only $\beta>2$. It would be possible also to improve the regularity condition on the initial data by working in Besov space invariant for the scaling of the system, but it is not the object of this paper. %As we are only interested by justifying where exactly in the proof of \cite{VG} we can improve the condition $\beta>2$, we are just going to give the proof of the estimate $L^{\infty}_{T}(L^{p}(\mathbb{T}^{2}))$ for the density with a such condition (see the p1115 in \cite{VG}).
\end{remarka}
We are going to give your second result on Korteweg system with supercritical smallness condition on the initial data; before let us give the following definition:
\begin{definition}
We set $q=\rho-1$, $m=\rho u$ and $\rho=h$.
\end{definition}
In the following we are dealing with the euclidian space $\R^{N}$ with $N\geq 2$. Let us give our main result on the Korteweg system where we prove the existence of global strong solution with large initial data on the irrotational part.
%\begin{theorem}
%If $u_{0}\in L^{2+\e}(\R^{2})$ and $\rho_{0}\in L^{\infty}(\R^{2})$ then the system (\ref{1}) has at least one global weak solution $(\rho,u)$ which verifies the following estimates
%\end{theorem}
\begin{theorem}
Suppose that we are under the conditions (\ref{hypothese}). Assume that $m_{0}\in B^{\N-1}_{2,1}$ and $q_{0}\in B^{\N}_{2,1}$ with $h_{0}\geq c>0$.
Then there exists a constant $\e_{0}$ depending on $\frac{1}{h_{0}}$ such that if:
$$\|m_{0}\|_{B^{\N-1}_{2,1}}\leq\e_{0},$$
then there exists a unique global solution $(q,m)$ for system (\ref{3system3})
with $h$ bounded away from zero and,
$$
\begin{aligned}
&h\in \widetilde{C}(\R^{+},B^{\N}_{2,1}%\cap B^{\frac{N}{p_{2}}-1}_{p_{1},1}
)\cap L^{1}(\R^{+},
B^{\N+2}_{2,1})\;\;\;\mbox{and}\;\;\;\;m\in \widetilde{C}(\R^{+};B^{\N-1}_{2,1})
\cap L^{1}(\R{+},B^{\N-1}_{2,1}\cap B^{\N+1}_{2,1}).
\end{aligned}
$$
\label{theo2}
\end{theorem}
%In the following theorem we are going to assume that we have no friction term.
%\begin{theorem}
%Suppose that we are under the conditions (\ref{hypothese}). Assume that $m_{0}\in B^{\N-1}_{2,1}$ and $q_{0}\in \widetilde{B}^{\N-1,\N}_{2,1}$ with $\rho_{0}\geq c>0$.
%Then there exists a constant $\e_{0}$ depending on $\frac{1}{h_{0}}$ such that if:
%$$\|m_{0}\|_{B^{\N-1}_{2,1}}\leq\e_{0},$$
%then there exists a unique global solution $(q,m)$ for system (\ref{3system3})
%with $h$ bounded away from zero and,
%$$
%\begin{aligned}
%&h\in \widetilde{C}(\R^{+},B^{\N}_{2,1}%\cap B^{\frac{N}{p_{2}}-1}_{p_{1},1}
%)\cap L^{1}(\R^{+},
%B^{\N+2}_{2,1})\;\;\;\mbox{and}\;\;\;\;m\in \widetilde{C}(\R^{+};B^{\N-1}_{2,1})
%\cap L^{1}(\R{+},B^{\N-1}_{2,1}\cap B^{\N+1}_{2,1}).
%\end{aligned}
%$$
%\label{theo3}
%\end{theorem}
%\begin{remarka}
%This last method has been employed in \cite{}
%\end{remarka}
Let us give the plane of this paper, we shall remind in section \ref{section1} some auxiliary results of Gagliardo-Nirenberg's inequality and in section \ref{section2} the Litllewood-Paley theory. In section \ref{section4} and section \ref{section5}, we will prove different a priori estimates on the density and the velocity which show the theorem \ref{theo1}. We will conclude in section \ref{section6} by the proof of theorem  \ref{theo2}.
\subsubsection*{Notation}
In all the paper, $C$ will stand for a ÒharmlessÓ constant, and we will sometimes use the notation $A \lesssim B$ equivalently to $A\leq CB$.
\section{Auxiliary Assertions}
\label{section1}
We are going to recall some lemma which are also present in \cite{VG} and that we prefer to state for the sake of the completnes.
\begin{lemme}
\label{lemme1}
 Let $\Omega\in\R^{N}$ be an arbitrary bounded domain satisfying the cone condition. Then the following inequality is valid for every function $u\in W^{1,m}(\Omega)$, $\int_{\Omega}u dx=0$
\begin{equation}
\|u\|_{L^{q}(\Omega)}\leq C_{1}\|\n u\|_{L^{m}(\Omega)}^{\alpha}\|u\|^{1-\alpha}_{L^{r}(\Omega)},
 \label{6}
\end{equation}
where $\alpha=\frac{\frac{1}{r}-\frac{1}{q}}{\frac{1}{r}-\frac{1}{m}+\frac{1}{n}}$, moreover if $m<n$ then $q\in[r,\frac{mn}{n-m}]$ for $r\leq\frac{mn}{n-m}$ and $q\in[\frac{mn}{n-m},r]$ for $r>\frac{mn}{n-m}$. If $m\geq n$ then $q\in[r,+\infty0$ is arbitrary; moreover if $m>n$ then equality (\ref{6}) is also valid for $q=+\infty$.
\end{lemme}
Inequality (\ref{6}) is a particular case of the more general inequalities proven in \cite{20,21,22}. %In \cite{23}, inequality (\ref{6}) is proven for an arbitrary $m\geq1$, moreover the dependence of the constant $C_{1}$ is indicated on the parameter involved in (\ref{6}).\\
Let us mention that an inequality of the form (\ref{6}) is valid for the function of class $W^{1,m}(\Omega)$ when $M=\frac{1}{|\Omega|}\int_{\Omega}u\,dx$ is not null. It suffices to consider $v=u-M$ and apply inequality (\ref{6}) to the function $v$. We obtain then the inequality
\begin{equation}
\|u\|_{L^{q}(\Omega)}\leq C_{2}(\|\n u\|_{L^{m}(\Omega)}^{\alpha}\|u\|^{1-\alpha}_{L^{r}(\Omega)}+\|u\|_{L^{1}(\Omega)}),
 \label{7}
\end{equation}
\begin{lemme}
Let $\Omega\in\R^{2}$ be an arbitrary bounded domain satisfying the cone condition. Then every function $u\in W^{1,m}(\Omega)$ with $\int_{\Omega}udx=0$ satisfies the inequality
\begin{equation}
\|u\|_{L^{\frac{2m}{2-m}}(\Omega)}\leq C_{3}(2-m)^{-\frac{1}{2}}\|\n u\|_{L^{m}(\Omega)},\;\;\;1\leq m<2,
 \label{8}
\end{equation}
where $C_{3}$ is a constant independent of $m$ and the function $u$.
 \label{lemme2}
\end{lemme}
For a proof of this inequality see \cite{24,26}. The exact constant in inequality (\ref{8}) is obtained in the article \cite{26}.
\begin{lemme}
Let $\Omega\in\R^{2}$ be an arbitrary bounded domain satisfying the cone condition. Then for an arbitrary number $\e$, $1\geq 2\e\geq 0$, every function $h\in W^{1,\frac{2m}{m+\delta}}(\Omega)$, $m\geq 2$, $1\geq\delta\geq 0$, satisfies the inequality
\begin{equation}
\|h\|_{L^{2m}(\Omega)}\leq C_{4}(\|h\|_{L^{1}(\Omega)}+m^{\frac{1}{2}}\|\n h\|^{1-s}_{L^{\frac{2m}{m+\delta}}(\Omega)}
\|h\|^{s}_{L^{2(1-\e)}(\Omega)}),
\label{9}
\end{equation}
where $s=(1-\e)\frac{1-\delta}{m-\delta(1-\e)}$ and  $C_{4}$ is a positive constant independent of $m,\e,\delta$ and the function $h$.
 \label{lemme3}
\end{lemme}
%In the sequel, we will need from the following Gr\"onwall lemma (see \cite{VG} for the proof).
%\begin{lemme}
%Let $g(t)$ be a nonnegative function of class $L^{1}(0,T)$ and let $a(\e)$ and $b(\e)$ be nonnegative functions given for $\e>0$ and satisfying the condition $\e b(\e)a^{\e}(\e)\rightarrow 0$ as $\e\rightarrow0$. If a continuous function $Y(t)$ on $[0,T]$
%satisfies the inequality
%$$Y(t)\leq a(\e)+b(\e)\int^{t}_{0}g(\tau)Y^{1+\e}(\tau)d\tau,$$
%for every $\e>0$ then
%$$Y(t)\leq 2^{\frac{1}{\e_{0}}}a(\e_{0}),\;\;t\in[0,T],$$
%where $\e_{0}$ is a positive number satisfying the condition
%$$2\e_{0}b(\e_{0})a^{\e_{0}}(\e_{0})\int^{T}_{0}g(t)dt\leq 1.$$
% \label{lemme4}
%\end{lemme}
%\subsection{Besov space}
\section{Littlewood-Paley theory and Besov spaces}
\label{section2}
Throughout the paper, $C$ stands for a constant whose exact meaning depends on the context. The notation $A\lesssim B$ means
that $A\leq CB$.
For all Banach space $X$, we denote by $C([0,T],X)$ the set of continuous functions on $[0,T]$ with values in $X$.
For $p\in[1,+\infty]$, the notation $L^{p}(0,T,X)$ or $L^{p}_{T}(X)$ stands for the set of measurable functions on $(0,T)$
with values in $X$ such that $t\rightarrow\|f(t)\|_{X}$ belongs to $L^{p}(0,T)$.
Littlewood-Paley decomposition  corresponds to a dyadic
decomposition  of the space in Fourier variables.
We can use for instance any $\varphi\in C^{\infty}(\R^{N})$,
supported in
${\cal{C}}=\{\xi\in\R^{N}/\frac{3}{4}\leq|\xi|\leq\frac{8}{3}\}$
such that:
$$\sum_{l\in\mathbb{Z}}\varphi(2^{-l}\xi)=1\,\,\,\,\mbox{if}\,\,\,\,\xi\ne 0.$$
Denoting $h={\cal{F}}^{-1}\varphi$, we then define the dyadic
blocks by:
$$\D_{l}u=\varphi(2^{-l}D)u=2^{lN}\int_{\R^{N}}h(2^{l}y)u(x-y)dy\,\,\,\,\mbox{and}\,\,\,S_{l}u=\sum_{k\leq
l-1}\D_{k}u\,.$$ Formally, one can write that:
$$u=\sum_{k\in\mathbb{Z}}\D_{k}u\,.$$
This decomposition is called homogeneous Littlewood-Paley
decomposition. Let us observe that the above formal equality does
not hold in ${\cal{S}}^{'}(\R^{N})$ for two reasons:
\begin{enumerate}
\item The right hand-side does not necessarily converge in
${\cal{S}}^{'}(\R^{N})$.
\item Even if it does, the equality is not
always true in ${\cal{S}}^{'}(\R^{N})$ (consider the case of the
polynomials).
\end{enumerate}
\subsection{Homogeneous Besov spaces and first properties}
\begin{definition}
For
$s\in\R,\,\,p\in[1,+\infty],\,\,q\in[1,+\infty],\,\,\mbox{and}\,\,u\in{\cal{S}}^{'}(\R^{N})$
we set:
$$\|u\|_{B^{s}_{p,q}}=(\sum_{l\in\mathbb{Z}}(2^{ls}\|\D_{l}u\|_{L^{p}})^{q})^{\frac{1}{q}}.$$
The Besov space $B^{s}_{p,q}$ is the set of temperate distribution $u$ such that $\|u\|_{B^{s}_{p,q}}<+\infty$.
\end{definition}
\begin{remarka}The above definition is a natural generalization of the
nonhomogeneous Sobolev and H$\ddot{\mbox{o}}$lder spaces: one can show
that $B^{s}_{\infty,\infty}$ is the nonhomogeneous
H$\ddot{\mbox{o}}$lder space $C^{s}$ and that $B^{s}_{2,2}$ is
the nonhomogeneous space $H^{s}$.
\end{remarka}
\begin{proposition}
\label{derivation,interpolation}
The following properties holds:
\begin{enumerate}
\item there exists a constant universal $C$
such that:\\
$C^{-1}\|u\|_{B^{s}_{p,r}}\leq\|\n u\|_{B^{s-1}_{p,r}}\leq
C\|u\|_{B^{s}_{p,r}}.$
\item If
$p_{1}<p_{2}$ and $r_{1}\leq r_{2}$ then $B^{s}_{p_{1},r_{1}}\hookrightarrow
B^{s-N(1/p_{1}-1/p_{2})}_{p_{2},r_{2}}$.
\item $B^{s^{'}}_{p,r_{1}}\hookrightarrow B^{s}_{p,r}$ if $s^{'}> s$ or if $s=s^{'}$ and $r_{1}\leq r$.
\end{enumerate}
\label{interpolation}
\end{proposition}
Let now recall a few product laws in Besov spaces coming directly from the paradifferential calculus of J-M. Bony
(see \cite{BJM}) and rewrite on a generalized form in \cite{3BC}.
\begin{proposition}
\label{produit1}
We have the following laws of product:
\begin{itemize}
\item For all $s\in\R$, $(p,r)\in[1,+\infty]^{2}$ we have:
\begin{equation}
\|uv\|_{B^{s}_{p,r}}\leq
C(\|u\|_{L^{\infty}}\|v\|_{B^{s}_{p,r}}+\|v\|_{L^{\infty}}\|u\|_{B^{s}_{p,r}})\,.
\label{2.2}
\end{equation}
\item Let $(p,p_{1},p_{2},r,\lambda_{1},\lambda_{2})\in[1,+\infty]^{2}$ such that:$\frac{1}{p}\leq\frac{1}{p_{1}}+\frac{1}{p_{2}}$,
$p_{1}\leq\lambda_{2}$, $p_{2}\leq\lambda_{1}$, $\frac{1}{p}\leq\frac{1}{p_{1}}+\frac{1}{\lambda_{1}}$ and
$\frac{1}{p}\leq\frac{1}{p_{2}}+\frac{1}{\lambda_{2}}$. We have then the following inequalities:\\
if $s_{1}+s_{2}+N\inf(0,1-\frac{1}{p_{1}}-\frac{1}{p_{2}})>0$, $s_{1}+\frac{N}{\lambda_{2}}<\frac{N}{p_{1}}$ and
$s_{2}+\frac{N}{\lambda_{1}}<\frac{N}{p_{2}}$ then:
\begin{equation}
\|uv\|_{B^{s_{1}+s_{2}-N(\frac{1}{p_{1}}+\frac{1}{p_{2}}-\frac{1}{p})}_{p,r}}\lesssim\|u\|_{B^{s_{1}}_{p_{1},r}}
\|v\|_{B^{s_{2}}_{p_{2},\infty}},
\label{2.3}
\end{equation}
when $s_{1}+\frac{N}{\lambda_{2}}=\frac{N}{p_{1}}$ (resp $s_{2}+\frac{N}{\lambda_{1}}=\frac{N}{p_{2}}$) we replace
$\|u\|_{B^{s_{1}}_{p_{1},r}}\|v\|_{B^{s_{2}}_{p_{2},\infty}}$ (resp $\|v\|_{B^{s_{2}}_{p_{2},\infty}}$) by
$\|u\|_{B^{s_{1}}_{p_{1},1}}\|v\|_{B^{s_{2}}_{p_{2},r}}$ (resp $\|v\|_{B^{s_{2}}_{p_{2},\infty}\cap L^{\infty}}$),
if $s_{1}+\frac{N}{\lambda_{2}}=\frac{N}{p_{1}}$ and $s_{2}+\frac{N}{\lambda_{1}}=\frac{N}{p_{2}}$ we take $r=1$.
\\
If $s_{1}+s_{2}=0$, $s_{1}\in(\frac{N}{\lambda_{1}}-\frac{N}{p_{2}},\frac{N}{p_{1}}-\frac{N}{\lambda_{2}}]$ and
$\frac{1}{p_{1}}+\frac{1}{p_{2}}\leq 1$ then:
\begin{equation}
\|uv\|_{B^{-N(\frac{1}{p_{1}}+\frac{1}{p_{2}}-\frac{1}{p})}_{p,\infty}}\lesssim\|u\|_{B^{s_{1}}_{p_{1},1}}
\|v\|_{B^{s_{2}}_{p_{2},\infty}}.
\label{2.4}
\end{equation}
If $|s|<\NN$ for $p\geq2$ and $-\frac{N}{p^{'}}<s<\NN$ else, we have:
\begin{equation}
\|uv\|_{B^{s}_{p,r}}\leq C\|u\|_{B^{s}_{p,r}}\|v\|_{B^{\NN}_{p,\infty}\cap L^{\infty}}.
\label{2.5}
\end{equation}
\end{itemize}
\end{proposition}
\begin{remarka}
In the sequel $p$ will be either $p_{1}$ or $p_{2}$ and in this case $\frac{1}{\lambda}=\frac{1}{p_{1}}-\frac{1}{p_{2}}$
if $p_{1}\leq p_{2}$, resp $\frac{1}{\lambda}=\frac{1}{p_{2}}-\frac{1}{p_{1}}$
if $p_{2}\leq p_{1}$.
\end{remarka}
\begin{corollaire}
\label{produit2}
Let $r\in [1,+\infty]$, $1\leq p\leq p_{1}\leq +\infty$ and $s$ such that:
\begin{itemize}
\item $s\in(-\frac{N}{p_{1}},\frac{N}{p_{1}})$ if $\frac{1}{p}+\frac{1}{p_{1}}\leq 1$,
\item $s\in(-\frac{N}{p_{1}}+N(\frac{1}{p}+\frac{1}{p_{1}}-1),\frac{N}{p_{1}})$ if $\frac{1}{p}+\frac{1}{p_{1}}> 1$,
\end{itemize}
then we have if $u\in B^{s}_{p,r}$ and $v\in B^{\frac{N}{p_{1}}}_{p_{1},\infty}\cap L^{\infty}$:
$$\|uv\|_{B^{s}_{p,r}}\leq C\|u\|_{B^{s}_{p,r}}\|v\|_{B^{\frac{N}{p_{1}}}_{p_{1},\infty}\cap L^{\infty}}.$$
\end{corollaire}
The study of non stationary PDE's requires space of type $L^{\rho}(0,T,X)$ for appropriate Banach spaces $X$. In our case, we
expect $X$ to be a Besov space, so that it is natural to localize the equation through Littlewood-Payley decomposition. But, in doing so, we obtain
bounds in spaces which are not type $L^{\rho}(0,T,X)$ (except if $r=p$).
We are now going to
define the spaces of Chemin-Lerner in which we will work, which are
a refinement of the spaces
$L_{T}^{\rho}(B^{s}_{p,r})$.
$\hspace{15cm}$
\begin{definition}
Let $\rho\in[1,+\infty]$, $T\in[1,+\infty]$ and $s_{1}\in\R$. We set:
$$\|u\|_{\widetilde{L}^{\rho}_{T}(B^{s_{1}}_{p,r})}=
\big(\sum_{l\in\mathbb{Z}}2^{lrs_{1}}\|\D_{l}u(t)\|_{L^{\rho}(L^{p})}^{r}\big)^{\frac{1}{r}}\,.$$
We then define the space $\widetilde{L}^{\rho}_{T}(B^{s_{1}}_{p,r})$ as the set of temperate distribution $u$ over
$(0,T)\times\R^{N}$ such that %$\lim_{q\rightarrow+\infty}S_{q}u=0$ in ${\cal S}^{'}((0,T)\times\R^{N})$
%and
$\|u\|_{\widetilde{L}^{\rho}_{T}(B^{s_{1}}_{p,r})}<+\infty$.
\end{definition}
We set $\widetilde{C}_{T}(\widetilde{B}^{s_{1}}_{p,r})=\widetilde{L}^{\infty}_{T}(\widetilde{B}^{s_{1}}_{p,r})\cap
{\cal C}([0,T],B^{s_{1}}_{p,r})$.
Let us emphasize that, according to Minkowski inequality, we have:
$$\|u\|_{\widetilde{L}^{\rho}_{T}(B^{s_{1}}_{p,r})}\leq\|u\|_{L^{\rho}_{T}(B^{s_{1}}_{p,r})}\;\;\mbox{if}\;\;r\geq\rho
,\;\;\;\|u\|_{\widetilde{L}^{\rho}_{T}(B^{s_{1}}_{p,r})}\geq\|u\|_{L^{\rho}_{T}(B^{s_{1}}_{p,r})}\;\;\mbox{if}\;\;r\leq\rho
.$$
\begin{remarka}
It is easy to generalize proposition \ref{produit1},
to $\widetilde{L}^{\rho}_{T}(B^{s_{1}}_{p,r})$ spaces. The indices $s_{1}$, $p$, $r$
behave just as in the stationary case whereas the time exponent $\rho$ behaves according to H\"older inequality.
\end{remarka}
In the sequel we will need of composition lemma in $\widetilde{L}^{\rho}_{T}(B^{s}_{p,r})$ spaces.
\begin{lemme}
\label{composition}
Let $s>0$, $(p,r)\in[1,+\infty]$ and $u\in \widetilde{L}^{\rho}_{T}(B^{s}_{p,r})\cap L^{\infty}_{T}(L^{\infty})$.
\begin{enumerate}
 \item Let $F\in W_{loc}^{[s]+2,\infty}(\R^{N})$ such that $F(0)=0$. Then $F(u)\in \widetilde{L}^{\rho}_{T}(B^{s}_{p,r})$. More precisely there exists a function $C$ depending only on $s$, $p$, $r$, $N$ and $F$ such that:
$$\|F(u)\|_{\widetilde{L}^{\rho}_{T}(B^{s}_{p,r})}\leq C(\|u\|_{L^{\infty}_{T}(L^{\infty})})\|u\|_{\widetilde{L}^{\rho}_{T}(B^{s}_{p,r})}.$$
\item If $v,\,u\in\widetilde{L}^{\rho}_{T}(B^{s}_{p,r})\cap
L^{\infty}_{T}(L^{\infty})$ and $G\in
W^{[s]+3,\infty}_{loc}(\R^{N})$ then $G(u)-G(v)$ belongs to
$\widetilde{L}^{\rho}_{T}(B^{s}_{p})$ and there exists a constant C
depending only of $s, p ,N\;\mbox{and}\;G$ such that:
$$
\begin{aligned}
\|G(u)-G(v)\|_{\widetilde{L}^{\rho}_{T}(B^{s}_{p,r})}\leq&
\,\,C(\|u\|_{L^{\infty}_{T}(L^{\infty})},\|v\|_{L^{\infty}_{T}(L^{\infty})})
\big(\|v-u\|_{\widetilde{L}^{\rho}_{T}(B^{s}_{p,r})}
(1+\|u\|_{L^{\infty}_{T}(L^{\infty})}\\
&+\|v\|_{L^{\infty}_{T}(L^{\infty})})+\|v-u\|_{L^{\infty}_{T}(L^{\infty})}(\|u\|_{\widetilde{L}^{\rho}_{T}(B^{s}_{p,r})}
+\|v\|_{\widetilde{L}^{\rho}_{T}(B^{s}_{p,r})})\big).
\end{aligned}
$$
%Let $F\in W_{loc}^{[s]+3,\infty}(\R^{N})$ such that $F(0)=0$. Then $F(u)-F^{'}(0)u\in \widetilde{L}^{\rho}_{T}(B^{s}_{p,r})$. More precisely there exists a function $C$ depending only on $s$, $p$, $r$, $N$ and $F$ such that:
%$$\|F(u)-F^{'}(0)u\|_{\widetilde{L}^{\rho}_{T}(B^{s}_{p,r})}\leq C(\|u\|_{L^{\infty}_{T}(L^{\infty})})\|u\|^{2}_{\widetilde{L}^{\rho}_{T}(B^{s}_{p,r})}.$$
\end{enumerate}
\end{lemme}
Now we give some result on the behavior of the Besov spaces via some pseudodifferential operator (see \cite{BCD}).
\begin{definition}
Let $m\in\R$. A smooth function function $f:\R^{N}\rightarrow\R$ is said to be a ${\cal S}^{m}$ multiplier if for all muti-index $\alpha$, there exists a constant $C_{\alpha}$ such that:
$$\forall\xi\in\R^{N},\;\;|\p^{\alpha}f(\xi)|\leq C_{\alpha}(1+|\xi|)^{m-|\alpha|}.$$
\label{smoothf}
\end{definition}
\begin{proposition}
Let $m\in\R$ and $f$ be a ${\cal S}^{m}$ multiplier. Then for all $s\in\R$ and $1\leq p,r\leq+\infty$ the operator $f(D)$ is continuous from $B^{s}_{p,r}$ to $B^{s-m}_{p,r}$.
\label{singuliere}
\end{proposition}
%Actually, in \cite{BCD}, the proposition below is proved for non-homogeneous Besov spaces. The adaptation to homogeneous spaces is straightforward.
Let us now give some estimates for the heat equation:
\begin{proposition}
\label{5chaleur} Let $s\in\R$, $(p,r)\in[1,+\infty]^{2}$ and
$1\leq\rho_{2}\leq\rho_{1}\leq+\infty$. Assume that $u_{0}\in B^{s}_{p,r}$ and $f\in\widetilde{L}^{\rho_{2}}_{T}
(B^{s-2+2/\rho_{2}}_{p,r})$.
Let u be a solution of:
$$
\begin{cases}
\begin{aligned}
&\p_{t}u-\mu\D u=f\\
&u_{t=0}=u_{0}\,.
\end{aligned}
\end{cases}
$$
Then there exists $C>0$ depending only on $N,\mu,\rho_{1}$ and
$\rho_{2}$ such that:
$$\|u\|_{\widetilde{L}^{\rho_{1}}_{T}(\widetilde{B}^{s+2/\rho_{1}}_{p,r})}\leq C\big(
 \|u_{0}\|_{B^{s}_{p,r}}+\mu^{\frac{1}{\rho_{2}}-1}\|f\|_{\widetilde{L}^{\rho_{2}}_{T}
 (B^{s-2+2/\rho_{2}}_{p,r})}\big)\,.$$
 If in addition $r$ is finite then $u$ belongs to $C([0,T],B^{s}_{p,r})$.
\end{proposition}
\section{Proof of theorem \ref{theo1}}
\label{section4}
In the sequel we shall work on the torus $\Omega=\mathbb{T}^{2}$. Let us start with recalling the energy estimate, when we multiply the momentum equation we get:
\begin{equation}
\begin{aligned}
&\int_{\Omega}(\rho|u|^{2}(t,x)+\pi(\rho)(t,x))dx+\int^{t}_{0}\int_{\Omega}|\n u|^{2}(s,x)dsdx\\
&+\int^{t}_{0}\int_{\Omega}(1+\lambda(\rho)(s,x))({\rm div}u)^{2}(s,x)dsdx\leq\int_{\Omega}(\rho_{0}(x)|u_{0}(x)|^{2}+\Pi(\rho_{0})(x))dx
\end{aligned}
\label{energy}
\end{equation}
with $\pi$ defined as follows:
$$
\pi(\rho)=a(\frac{1}{\gamma-1}(\rho^{\gamma}-\rho)-\rho+1)\:\;\mbox{for}\; \gamma>1.$$% \left\{ 
%\begin{aligned}
%&a(\frac{1}{\gamma-1}(\rho^{\gamma}-\rho)-\rho+1)\:\;\mbox{for}\; \gamma>1% 1,\;\gamma\ne 1\\
%&a(\rho\ln\rho-\rho+1)\:\;\mbox{for}\; \gamma=1.
%\end{aligned}$$
% \right\}. %\quad 1 \leq i \leq N. \hspace{8mm} (\rm{I}. 
%1)
%$$
Let us recall that $P^{'}(\rho)=\rho \pi^{''}(\rho)$ what implies by convexity that $\pi(\rho)\geq 0$. Finally as we assume that:
$$C_{1}=\int_{\Omega}\big(\frac{1}{2}\rho_{0}(x)|u_{0}(x)|^{2}+\pi(\rho_{0})(x)+\rho_{0}(x)\big)dx$$
is finite, we obtain at least formally (if $\rho$ and $u$ are enough regular for performing integration by parts) by energy estimate (\ref{energy}) and via the transport equation that:
\begin{equation}
\begin{aligned}
&\int_{\Omega}(\rho|u|^{2}(t,x)+\pi(\rho)(t,x)+\rho(t,x))dx+\int^{t}_{0}\int_{\Omega}|\n u|^{2}(s,x)dsdx\\
&\hspace{5cm}+\int^{t}_{0}\int_{\Omega}(1+\lambda(\rho)(s,x))({\rm div}u)^{2}(s,x)dsdx\leq C_{1}
\end{aligned}
\label{energy1}
\end{equation}
Let us now explain how to get $L^{2}((0,T)\times\Omega)$ estimates on $u$, we are going to follow Lions in \cite{PL2} p4. Indeed by the momentum equation we have:
$$|\int_{\Omega}\rho u(t,x)dx|=|\int_{\Omega}\rho_{0} u_{0}(x)dx|\leq\|\rho_{0}u_{0}\|_{L^{1}(\Omega)}.$$
Next we use the Poincar\'e-Wirtinger inequality and we have:
$$|\int_{\Omega}\rho(t,x)[u(t,x)-\int_{\Omega}u(t,y)dy]dx|\leq C\|\rho(t,\cdot)\|_{L^{\gamma}}\|\n u\|_{L^{2}(\Omega)}.$$
Hence for all $t\geq 0$:
$$|\int_{\Omega}u(t,x)dx|\leq\frac{1}{(\int_{\Omega}\rho_{0}dx)}(\|\rho_{0}u_{0}\|_{L^{1}(\Omega)}+
C\|\rho(t,\cdot)\|_{L^{\gamma}}\|\n u\|_{L^{2}(\Omega)}).$$
We conclude by Poincar\'e-Wirtinger inequality which implies that $|\int_{\Omega}u(t,x)dx|+\|\n u\|_{L^{2}(\Omega)}$ is an equivalent norm to the usual one in $H^{1}(\Omega)$.
\hfill {$\Box$}\\
\\
Now we are just going to explain where in the proof we can slightly improve the range of the coefficient $\beta$ in \cite{VG}. One of the main point of the proof in \cite{VG} consists in getting a priori estimates on the density in $L^{\infty}(L^{p}(\T^{2}))$ for any $p>1$. This is possible due to the viscosity coefficient $\lambda(\rho)=\rho^{\beta}$ which provide such estimate at least if $\beta$ is large enough. Let us follow the arguments of the proof of \cite{VG} and explain where by commutators estimates we can weaken the hypothesis $\beta>3$.
\section{A priori estimates on the density and the velocity}
\label{section5}
First as in \cite{VG}, we are going to recall some estimates for solutions to the following two Neumann problems:
\begin{equation}
\D\xi={\rm div}(\rho u),\;\;\int_{\Omega}\xi dx=0,\;\;\p_{x_{1}}\xi |_{x_{1}=0,x_{1}=1}=\p_{x_{2}}\xi |_{x_{2}=0,x_{2}=1}=0.
\label{19}
\end{equation}
\begin{equation}
\D\eta={\rm div}({\rm div}(\rho u\otimes u)),\;\;\int_{\Omega}\eta dx=0,,\;\;\p_{x_{1}}\eta |_{x_{1}=0,x_{1}=1}=\p_{x_{2}}\eta |_{x_{2}=0,x_{2}=1}=0.
\label{20}
\end{equation}
Therefore by \cite{27} we have solution to the problems (\ref{19}) and (\ref{20}), whereas the estimates for singular integrals in \cite{28} provide the following inequalities:
$$
\begin{aligned}
&\|\n(\D)^{-1}{\rm div}(\rho u)\|_{L^{2m}}\lesssim m\|\rho u\|_{L^{2m}},\;\;1\leq m<+\infty,\\
&\|\n(\D)^{-1}{\rm div}(\rho u)\|_{L^{2-r}}\lesssim \|\rho u\|_{L^{2-r}},\;\;1\geq 2r\geq 0,\\
&\|R_{i,j}(\rho u_{i}u_{j})\|_{L^{2m}}\lesssim m\|\rho u\otimes u\|_{L^{2m}},\;\;1\leq m<+\infty,\\
\end{aligned}
$$
Here we have roughly written $\xi=(\D)^{-1}{\rm div}(\rho u)$ and $\eta=R_{i,j}(\rho u_{i}u_{j})$ (with the summation notation).\\
By H\"older's inequalities we obtain:
%Making use of these inequalities together with $\rho>0$, we find the following estimates for all $m\geq1$ and $k>1$:
\begin{equation}
\begin{aligned}
&\|\n(\D)^{-1}{\rm div}(\rho u)\|_{L^{2m}}\lesssim m\|\rho\|_{L^{\frac{2mk}{k-1}}}\|u\|_{L^{2mk}},\\
&\|\n(\D)^{-1}{\rm div}(\rho u)\|_{L^{2-r}}\lesssim \|\rho\|_{L^{\frac{2-r}{r}}}\|\sqrt{\rho}u\|^{\frac{1}{2}}_{L^{2}},\\
&\|R_{i,j}(\rho u_{i}u_{j})\|_{L^{2m}}\lesssim m\|\rho\|_{L^{\frac{2mk}{k-1}}}\|u\|^{2}_{L^{4mk}},\\
\end{aligned}
\label{21}
\end{equation}
where $k>1$, $m\geq1$ and $r\geq1$, $1\geq 2r\geq0$.\\
From the estimate of lemma \ref{lemme1}-\ref{lemme3}, we obtain:
\begin{equation}
\begin{aligned}
&\|u\|_{L^{2m}}\lesssim m^{\frac{1}{2}}\|\n u\|_{L^{2}},\;\;m>2,\\
&\|(\D)^{-1}{\rm div}(\rho u)\|_{L^{2m}}\lesssim m^{\frac{1}{2}}\|\n (\D)^{-1}{\rm div}(\rho u)\|_{\frac{2m}{m+1}},\;\;m>2,
\end{aligned}
\label{22}
\end{equation}
We set now:
$$\va(t)=\int_{\Omega}\big({\rm curl}u^{2}(t,x)+(2+\lambda(\rho)){\rm div}u^{2}(t,x)\big)dx,$$
we obtain then by using (\ref{22}), (\ref{21}) with $r=\frac{2}{m+1}$ and the energy inequality (\ref{energy1}):
\begin{equation}
\begin{aligned}
&\|(\D)^{-1}{\rm div}(\rho u)\|_{L^{2m}}\lesssim m^{\frac{1}{2}}\|\rho\|_{L^{m}}^{\frac{1}{2}},\;\;m>2,\\
%&\|\n(\D)^{-1}{\rm div}(\rho u)\|_{L^{2m}}\lesssim m^{\frac{3}{2}}k^{\frac{1}{2}}(\va(t))^{\frac{1}{2}}\|\rho\|_{L^{\frac{2mk}{k-1}}},\;\;m>2,\;\;k>1,\\
%&\|R_{i,j}(\rho u_{i}u_{j})\|_{L^{2m}}\lesssim m^{2}k\va(t)\|\rho\|_{L^{\frac{2mk}{k-1}}},\;\;m>2,\;\;k>1.\\
%&\|[R_{i,j},u_{j}](\rho u_{j})\|_{L^{2m}}\lesssim m^{2}k\va(t)\|\rho\|_{L^{\frac{2mk}{k-1}}},\;\;m>2,\;\;k>1.
\end{aligned}
\label{23.1}
\end{equation}
Similarly we have:
\begin{equation}
\begin{aligned}
%&\|(\D)^{-1}{\rm div}(\rho u)\|_{L^{2m}}\lesssim m^{\frac{1}{2}}\|\rho\|_{L^{m}}^{\frac{1}{2}},\;\;m>2,\\
&\|\n(\D)^{-1}{\rm div}(\rho u)\|_{L^{2m}}\lesssim m^{\frac{3}{2}}k^{\frac{1}{2}}(\va(t))^{\frac{1}{2}}\|\rho\|_{L^{\frac{2mk}{k-1}}},\;\;m>2,\;\;k>1,\\
&\|R_{i,j}(\rho u_{i}u_{j})\|_{L^{2m}}\lesssim m^{2}k\va(t)\|\rho\|_{L^{\frac{2mk}{k-1}}},\;\;m>2,\;\;k>1.\\
%&\|[R_{i,j},u_{j}](\rho u_{j})\|_{L^{2m}}\lesssim m^{2}k\va(t)\|\rho\|_{L^{\frac{2mk}{k-1}}},\;\;m>2,\;\;k>1.
\end{aligned}
\label{23.2}
\end{equation}
\subsection{Gain of integrability for the density}
Following \cite{VG} the plan of the proof of \cite{VG}, we are interested  in getting a gain of integrability on the density. We follow here the method of Lions in \cite{PL2} to get a gain of integrability on the pressure and the argument developed in \cite{VG}. Apply the operator $(\D)^{-1}{\rm div}$ to the momentum equation, we obtain:
\begin{equation}
\begin{aligned}
\frac{\p}{\p t}(\D)^{-1}{\rm div}(\rho u)&+[R_{ij},u_{j}](\rho u_{i})-(2+\lambda(\rho)){\rm div}u\\
&+P(\rho)-\frac{1}{|\Omega|}\int_{\Omega}(P(\rho)(t,x)-(2+\lambda(\rho)){\rm div}u)dx=0.
\end{aligned}
\label{26}
\end{equation}
We will set in the sequel:
\begin{equation}
B=(2+\lambda(\rho)){\rm div}u-P(\rho).
\label{pressionefficace}
\end{equation}
Next if we renormalize the mass equation we have:
$$\p_{t}\theta(\rho)+u\cdot\n\theta(\rho)+\rho\theta^{'}(\rho){\rm div}u=0.$$
where we have set:
$$\theta(\rho)=\int^{\rho}_{1}\frac{1}{s}(2+\lambda(s))ds=2\ln\rho+\frac{1}{\beta}(\rho^{\beta}-1).$$%\ln\rho 1_{\rho\leq m}+\frac{1}{\beta}\rho^{\beta}1_{\rho\geq 2m}+\int^{\rho}_{1}\frac{\va(s)1_{\{m\leq s\leq 2m\}}}{s}ds,$$
Finally we get the following transport equation:
\begin{equation}
\begin{aligned}
\frac{\p}{\p t}\big[(\D)^{-1}{\rm div}(\rho u)&+\theta(\rho)\big]+u\cdot\n\big[(\D)^{-1}{\rm div}(\rho u)+\theta(\rho)\big]+[R_{ij},u_{j}](\rho u_{i})\\
&+P(\rho)%-P(\bar{\rho})=0.%
-\frac{1}{|\Omega|}\int_{\Omega}\big[P(\rho)(t,x)-(2+\lambda(\rho)){\rm div}u\big]dx=0,
\end{aligned}
\label{26}
\end{equation}
%In the sequel we will note:
%$B=-\frac{1}{|\Omega|}\int_{\Omega}P(\rho)(t,x)dx.$\\
%This is that special equation from which we infer the second a priori estimate for the density. 
Denote by $f$ the function:
$$f(t,x)=\max(0,(\D)^{-1}{\rm div}(\rho u)+\theta(\rho))$$
and multiply the equation (\ref{26}) by the function $\rho f^{2m-1}$ with $m\in\mathbb{N}$ and $m\geq 4$ and integrate  over $\Omega$
%$\R^{N}$
, we obtain:
\begin{equation}
\begin{aligned}
&\frac{1}{2m}\frac{d}{dt}\int_{\Omega}\rho f^{2m}dx+\int_{\Omega}\rho P(\rho)f^{2m-1}dx
+\int_{\Omega}[R_{ij},u_{j}](\rho u_{i})\rho f^{2m-1}dx\\
&\hspace{8cm}+\int_{\Omega} B\,dx \int_{\Omega}\rho f^{2m-1}dx=0.%\\
%&\hspace{1cm}+\int_{\Omega}Bdx\int_{\Omega}\rho f^{2m-1}dx=0\\
\end{aligned}
\label{28}
\end{equation}
As in \cite{VG} we set:
\begin{equation}
Z(t)=\big(\int_{\Omega}\rho f^{2m}(t,x)\big)^{\frac{1}{2m}}
\label{29}
\end{equation}
Using %inequalities (\ref{18}), (\ref{22}) and (\ref{23}), the estimate $\rho\geq0$, condition (\ref{2}) on the parameter $\beta$ and
H\"older's inequality ($\frac{2m-1}{2m}+\frac{1}{2m}=1$, $\frac{\beta}{2m\beta+1}+\frac{1}{2m(2m\beta+1)}=\frac{1}{2m}$), we begin with estimating the term $|\int_{\Omega}[R_{ij},u_{j}](\rho u_{i})\rho f^{2m-1}dx$ in  (\ref{28}) as follows:
\begin{equation}
\begin{aligned}
\big|\int_{\Omega}[R_{ij},u_{j}](\rho u_{i})\rho f^{2m-1}dx|&\leq\int_{\Omega}|[R_{ij},u_{j}](\rho u_{i})|\rho^{\frac{1}{2m}}\rho^{\frac{2m-1}{2m}} f^{2m-1}dx\\
&\leq\|\,|[R_{ij},u_{j}](\rho u_{i})|\rho^{\frac{1}{2m}}\|_{L^{2m}(\Omega)}Z^{2m-1}(t)\\
&\leq\|\rho\|_{L^{2m\beta+1}(\Omega)}^{\frac{1}{2m}}\|[R_{ij},u_{j}](\rho u_{i})\|_{L^{2m+\frac{1}{\beta}}(\Omega)}
Z^{2m-1}(t)\\
%&\leq Z^{2m-1}(t)(1+\va(t)^{\frac{1}{2}}+\|\rho\|^{\frac{\beta}{2}}_{L^{2m\beta+1}(\Omega)}\va(t)^{\frac{1}{2}}).
\end{aligned}
\label{crucial}
\end{equation}
Next we recall a result of R. Coifman, P.-L. Lions, Y. Meyer and S. Semmes in \cite{4M}, which says that the following
map:
$$
\begin{aligned}
W^{1,r_{1}}(\mathbb{T}^{N})^{N}\times L^{r_{2}}(\mathbb{T}^{N})^{N}&\rightarrow W^{1,r_{3}}(\mathbb{T}^{N})^{N}\\
(a,b)&\rightarrow[a_{j}, R_{i,j}]b_{i}
\end{aligned}
$$
is continuous for any $N\geq2$ as soon as $\frac{1}{r_{3}}=\frac{1}{r_{1}}+\frac{1}{r_{2}}$. Hence we obtain that
$[R_{ij},u_{j}](\rho u_{i})$ belongs in $W^{1,p}$ (where $\frac{1}{p}=\frac{1}{2}+\frac{1}{2(m+\frac{1}{2\beta})k}+\frac{k-1}{2(m+\frac{1}{2\beta})k}$ with $k>1$ and $p=2-\frac{2}{m+1+\frac{1}{2\beta}}<2$) with the following inequality:
\begin{equation}
\begin{aligned}
\|[R_{ij},u_{j}](\rho u_{i})\|_{W^{1,p}(\Omega)}&\leq C\|\n u\|_{L^{2}(\Omega)}\|u\|_{L^{2(m+\frac{1}{2\beta})k}(\Omega)}\|\rho\|_{L^{\frac{2(m+\frac{1}{2\beta})k}{k-1}}(\Omega)},\\
&\leq C\|\n u\|_{L^{2}(\Omega)}\|u\|_{L^{2(m+\frac{1}{2\beta})k}(\Omega)}\|\rho\|_{L^{2m\beta+1}(\Omega)},
\end{aligned}
\label{a1}
\end{equation}
 where we have choose $k$ such that $\frac{2(m+\frac{1}{2\beta})k}{k-1}=2m\beta+1$, let $k=\frac{2m\beta+1}{2m(\beta-1)+1-\frac{1}{\beta}}$. We verifies that $\frac{1}{q}=\frac{1}{p}-\frac{1}{2}=\frac{1}{2m+\frac{1}{\beta}}$.
Next by using lemma \ref{lemme2} and (\ref{a1}) we get:
$$
\begin{aligned}
&\|[R_{ij},u_{j}](\rho u_{i})(t,\cdot)\|_{L^{2m+\frac{1}{\beta}}(\Omega)}\lesssim \big( m^{\frac{1}{2}}
\|\n u(t,\cdot)\|_{L^{2}(\Omega)}\|u(t,\cdot)\|_{L^{2(m+\frac{1}{2\beta})k}(\Omega)}\|\rho(t,\cdot)\|_{L^{2m\beta+1}(\Omega)}\\
&\hspace{10cm}+|\int_{\Omega}[R_{ij},u_{j}](\rho u_{i})(t,x)dx|\big).
\end{aligned}
$$
We can easily bound the last term on the right hand side by using the continuity of the Riez transform in $L^{p}(\omega)$ with $1<p<+\infty$:
$$|\int_{\Omega}[R_{ij},u_{j}](\rho u_{i})(t,x)dx|\lesssim \|\rho(t,\cdot)\|_{L^{\gamma}(\Omega)}\|u(t,\cdot)\|^{2}_{H^{1}(\Omega)}.$$
%$$
%\begin{aligned}
%\|[R_{ij},u_{j}](\rho u_{i})\|_{L^{2m+\frac{1}{\beta}}(\Omega)}&\leq C\big((m+\frac{1}{2\beta})k)^{\frac{1}{2}}\|\n %u\|_{L^{2}(\Omega)}^{2}\|\rho\|_{L^{2m\beta+1}(\Omega)}
%\end{aligned}
%$$
By (\ref{22}) and the previous inequalites we obtain finally:
\begin{equation}
\|[R_{ij},u_{j}](\rho u_{i})(t,\cdot)\|_{L^{2m+\frac{1}{\beta}}(\Omega)}\leq C m\|\n u\|_{L^{2}(\Omega)}^{2}\|\rho\|_{L^{2m\beta+1}(\Omega)}+\|\rho(t,\cdot)\|_{L^{\gamma}(\Omega)}\|u(t,\cdot)\|^{2}_{H^{1}(\Omega)}.
\label{crucial1}
\end{equation}
We have then from (\ref{crucial}) and (\ref{crucial1}):
$$
\begin{aligned}
\big|\int_{\Omega}[R_{ij},u_{j}](\rho u_{i})\rho f^{2m-1}dx|&\lesssim (m\|\rho\|_{L^{2m\beta+1}(\Omega)}^{1+\frac{1}{2m}}\va(t)+\|\rho(t,\cdot)\|_{L^{\gamma}(\Omega)}\|u(t,\cdot)\|^{2}_{H^{1}(\Omega)})
Z^{2m-1}(t).
%&\leq Z^{2m-1}(t)(1+\va(t)^{\frac{1}{2}}+\|\rho\|^{\frac{\beta}{2}}_{L^{2m\beta+1}(\Omega)}\va(t)^{\frac{1}{2}}).
\end{aligned}
$$
Next as in \cite{VG} we get:
$$
\begin{aligned}
\big|\int_{\Omega}Bdx\int_{\Omega}\rho f^{2m-1}dx\big|&\lesssim Z^{2m-1}(t)\|\rho\|_{L^{1}}^{\frac{1}{2m}}\int_{\Omega}\big((2+\lambda)|{\rm div}u|+P)dx\\
&\lesssim Z^{2m-1}(t)\big(1+(\va(t))^{\frac{1}{2}}(\int_{\Omega}(2+\lambda(\rho))dx)^{\frac{1}{2}}\big)\\
&\lesssim Z^{2m-1}(t)(1+\va(t)^{\frac{1}{2}}+\|\rho\|^{\frac{\beta}{2}}_{L^{2m\beta+1}(\Omega)}\va(t)^{\frac{1}{2}}).
\end{aligned}
$$
Collecting all the above inequalities, we obtain:
\begin{equation}
\begin{aligned}
&Z(t)\lesssim 1+\int_{0}^{t}\|\rho\|_{L^{\gamma}(\Omega)}(\tau)\|u\|^{2}_{H^{1}(\Omega)}(\tau)d\tau+\int_{0}^{t}m\va(\tau)\|\rho\|_{L^{2m\beta+1}(\Omega)}^{1+\frac{1}{2m}}(\tau)d\tau\\
&\hspace{8cm}+\int^{t}_{0}
\va(\tau)^{\frac{1}{2}}\|\rho\|_{L^{2m\beta+1}(\Omega)}^{\frac{\beta}{2}}(\tau)d\tau.
\end{aligned}
\label{30}
\end{equation}
As we have seen that $u$ belongs in $L^{2}((0,t),H^{1}(\Omega))$ we have:
$$\int_{0}^{t}\|\rho\|_{L^{\gamma}(\Omega)}(\tau)\|u\|^{2}_{H^{1}(\Omega)}(\tau)d\tau
\lesssim1.$$
We obtain then:
\begin{equation}
\begin{aligned}
&Z(t)\lesssim 1+\int_{0}^{t}m\va(\tau)\|\rho\|_{L^{2m\beta+1}(\Omega)}^{1+\frac{1}{2m}}(\tau)d\tau
+\int^{t}_{0}
\va(\tau)^{\frac{1}{2}}\|\rho\|_{L^{2m\beta+1}(\Omega)}^{\frac{\beta}{2}}(\tau)d\tau.
\end{aligned}
\label{301}
\end{equation}
Next following \cite{VG} we introduce the measurable sets:
$$\Omega_{1}(t)=\{x\in\Omega/\rho\geq 2m^{'}\}\;\;\;
\mbox{and}\;\;\;\Omega_{2}(t)=\{x\in\Omega_{1}(t)/\theta(\rho)+(\D)^{-1}{\rm div}(\rho u)>0\}.$$
We then have:
\begin{equation}
\|\rho\|_{L^{2m\beta+1}(\Omega)}^{\beta}\lesssim\big(\int_{\Omega_{1}(t)}\rho^{2m\beta+1}dx\big)^{\frac{\beta}{2m\beta+1}}
+1,
\label{31}
\end{equation}
Moreover by the definition of the function $\theta(\rho)$, we have:
\begin{equation}
\big(\int_{\Omega_{1}(t)}\rho^{2m\beta+1}dx\big)^{\frac{\beta}{2m\beta+1}}\lesssim
\big(\int_{\Omega_{1}(t)}\rho\theta(\rho)^{2m}dx\big)^{\frac{\beta}{2m\beta+1}}
\label{32}
\end{equation}
Using the fact that on $\Omega_{1}(t)\setminus \Omega_{2}(t)$ we have
$0\leq\theta(\rho)\leq|(\D)^{-1}{\rm div}(\rho u)|$, we derive the following estimate:
$$
\begin{aligned}
&\int_{\Omega_{1}(t)}\rho\theta(\rho)^{2m}dx=\int_{\Omega_{2}(t)}\rho(\theta(\rho)+(\D)^{-1}{\rm div}(\rho u)-
(\D)^{-1}{\rm div}(\rho u))^{2m}dx\\
&\hspace{9cm}+\int_{\Omega_{1}(t)\setminus \Omega_{2}(t)}\rho\theta(\rho)^{2m}dx\\
&\leq 2^{2m-1}\big(\int_{\Omega_{2}(t)}\rho f(\rho)^{2m}dx+\int_{\Omega_{2}(t)}\rho|(\D)^{-1}{\rm div}(\rho u)|^{2m}dx\big)\\
&\hspace{7,5cm}+\int_{\Omega_{1}(t)\setminus \Omega_{2}(t)}\rho |(\D)^{-1}{\rm div}(\rho u)|^{2m}dx\\
&\leq 2^{2m}(Z^{2m}(t)+\int_{\Omega}\rho |(\D)^{-1}{\rm div}(\rho u)|^{2m}dx)
\end{aligned}
$$
From estimates (\ref{31}) and (\ref{32}) we deduce:
$$\|\rho\|_{L^{2m\beta+1}(\Omega)}^{\beta}\leq C\big((m^{'})^{\beta}+Z(t)^{\frac{2m\beta}{2m\beta+1}}(t)+(\int_{\Omega}\rho |(\D)^{-1}{\rm div}(\rho u)|^{2m}dx)^{\frac{\beta}{2m\beta+1}}\big)$$
In view of estimate (\ref{23.1}), (\ref{23.2}), we have:
$$
\begin{aligned}
\int_{\Omega}\rho |(\D)^{-1}{\rm div}(\rho u)|^{2m}dx&\leq\|\rho\|_{L^{2m\beta+1}(\Omega)}\|(\D)^{-1}{\rm div}(\rho u)\|^{2m}_{L^{2m+\frac{1}{\beta}}(\Omega)}\\
&\leq \|\rho\|_{L^{2m\beta+1}(\Omega)}\big( (m+\frac{1}{2\beta})^{\frac{1}{2}}\|\rho\|_{L^{m+\frac{1}{\beta}}(\Omega)}^{\frac{1}{2}}\big)^{2m}\\
&\leq C^{m}m^{m}\|\rho\|_{L^{2m\beta+1}(\Omega)}^{m+1}.
\end{aligned}
$$
Finaly:
\begin{equation}
\|\rho\|_{L^{2m\beta+1}(\Omega)}^{\beta}\lesssim \big(Z(t)^{\frac{2m\beta}{2m\beta+1}}(t)+m^{\frac{1}{2}}
\|\rho\|_{L^{2m\beta+1}(\Omega)}^{\frac{\beta(m+1)}{2m\beta+1}}\big).
\label{33}
\end{equation}
By Young's inequality (with $q=\frac{2m\beta+1}{m+1}$ and $p=\frac{2\beta}{2\beta-1}+\frac{1}{m(2\beta-1)}$), we obtain:
$$
\|\rho\|_{L^{2m\beta+1}(\Omega)}^{\beta}\lesssim \big(Z(t)+\frac{1}{\e}m^{\frac{\beta}{2\beta-1}+\frac{1}{2m(2\beta-1)}}+
\e\|\rho\|_{L^{2m\beta+1}(\Omega)}^{\beta}\big).
$$
By bootstrap, we get:
\begin{equation}
\|\rho\|_{L^{2m\beta+1}(\Omega)}^{\beta}\lesssim Z(t)+m^{\frac{\beta}{2\beta-1}}.
\label{34}
\end{equation}
Therefore (\ref{30}) and (\ref{34})give the following inequality:
$$
\|\rho\|_{L^{2m\beta+1}(\Omega)}^{\beta}\lesssim\big(m^{\frac{\beta}{2\beta-1}}+\int_{0}^{t}m\va(\tau)\|\rho\|_{L^{2m\beta+1}(\Omega)}^{1+\frac{1}{2m}}(\tau)d\tau+\int^{t}_{0}
\va(\tau)^{\frac{1}{2}}\|\rho\|_{L^{2m\beta+1}(\Omega)}^{\frac{\beta}{2}}(\tau)d\tau\big)
$$
Next by Young's inequality, we have:
$$
\|\rho\|_{L^{2m\beta+1}(\Omega)}^{\beta}\lesssim\big(1+m^{\frac{\beta}{2\beta-1}}+\int_{0}^{t}m\va(\tau)\|\rho\|_{L^{2m\beta+1}(\Omega)}^{1+\frac{1}{2m}}(\tau)d\tau+\int^{t}_{0}
\|\rho\|_{L^{2m\beta+1}(\Omega)}^{\beta}(\tau)d\tau\big)
$$
%Attention!!!!!!!!!!!!!!Ici les constantes d\'ependent du temps, plus le temps est grand plus elles sont grandes.\\
%\\
Using the fact that $\va(t)\in L^{1}(0,T)$ and applying Gr\"onwall's inequality, wehave that:
%$$
%\|\rho\|_{L^{2m\beta+1}(\Omega)}^{\beta}\leq C\big(1+m^{\frac{\beta}{\beta-1}}+\int_{0}^{t}m^{2}\va(\tau)\|\rho\|_{L^{2m\beta+1}(\Omega)}^{1+\frac{1}{2m}}(\tau)d\tau\big)$$%+
%\frac{1}{\e}\int^{t}_{0}
%\va(\tau)d\tau+\e\int^{t}_{0}\|\rho\|_{L^{2m\beta+1}(\Omega)}^{\beta}(\tau)d\tau\big)
%$$
%and:
$$
\|\rho\|_{L^{2m\beta+1}(\Omega)}^{\beta}\leq C\big(1+m^{\frac{\beta}{2\beta-1}}+\int_{0}^{t}m\va(\tau)\|\rho\|_{L^{2m\beta+1}
(\Omega)}^{1+\frac{1}{2m}}(\tau)d\tau,
\big)
$$
where $C$ depends on $t$.
Denote:
\begin{equation}
y(t)=m^{-\frac{1}{\beta-1}}\|\rho\|_{L^{2m\beta+1}(\Omega)},\;\;\;t\in[0,T].
\label{35}
\end{equation}
Then:
$$y^{\beta}(t)m^{\frac{\beta}{\beta-1}}\leq C( 1+m^{\frac{\beta}{2\beta-1}}+
m\,m^{(1+\frac{1}{2m})\frac{1}{\beta-1}}\int_{0}^{t}\va(\tau)y^{1+\frac{1}{2m}}(\tau)d\tau
\big)$$
and we have:
$$m m^{(1+\frac{1}{2m})\frac{1}{\beta-1}}=m^{\frac{\beta}{\beta-1}+\frac{1}{2m(\beta-1)}}.$$
We have then:
$$y^{\beta}(t)\leq C\big(1+m^{\frac{-\beta^{2}}{(\beta-1)(2\beta-1)}}+
m^{\frac{1}{2m(\beta-1)}}\int_{0}^{t}\va(\tau)y^{1+\frac{1}{2m}}(\tau)d\tau
\big)$$
where $\frac{\beta}{2\beta-1}-\frac{\beta}{\beta-1}=\frac{-\beta^{2}}{(\beta-1)(2\beta-1)}<0$.\\
Recalling that $\beta>1$ and $\va(t)\in L^{1}(0,T)$ we find that for $m$ big enough:
$$y^{\beta}(t)\leq C\big(C_{1}+
\int_{0}^{t}\va(\tau)y^{\beta}(\tau)d\tau
\big)$$
%Attention!!!!!!!!!! ici la constante $C$ est independante de $m^{'}$ seulement pour des $m$ tres grand car il ecrasent $m^{'}$, sinon pour des $m$ de la taille de $m^{'}$ c'est different et en particulier si $m=m^{'}$ alors on a du
%$(m^{'})^{\frac{\beta(\beta-3)}{\beta-1}}$.\\
whence by Gr\"onwall inequality:
$$y(t)\leq C% \big(m^{\frac{-1}{\beta-1}}+1)
,\;\;\;t\in[0,T],$$
where $C$ depends on $t$.
We thus have:
$$\|\rho\|_{L^{2m\beta+1}}\leq C m^{\frac{1}{\beta-1}},\;\;\;t\in[0,T].$$
Hence the inequality:
\begin{equation}
\|\rho\|_{L^{k}(\Omega)}(t)\leq C k^{\frac{1}{\beta-1}},\;\;t\in[0,T].
\label{36}
\end{equation}
is valid for every $k\geq1$, with $C$ a positive constant independent of $k\geq1$ but depending of the time. \hfill {$\Box$}
\begin{remarka}
Let us point out that the estimate (\ref{36}) is the key point in order to improve the range on $\beta$. Indeed this last one is a refinement of the corresponding one in \cite{VG}. In particular we will be able to obtain the energy estimates (\ref{57}) only with assuming $\beta>2$.
\end{remarka}
\subsection{Second a priori estimate for the velocity}
In this section, we are going to furnish estimates on the velocity by using the gain of integrability on the density proved
in the previous section. We are going essentially to follow the proof of \cite{VG} and to emphasize
on the key point where we will only need the hypothesis $\beta>2$. We begin with recalling some equation on the effective pressure defined in \cite{PL2} and the rotational $curl$. We set:
$$
\begin{aligned}
&A={\rm curl}u\;\;\;\mbox{and}\;\;\;B=(2+\lambda(\rho)){\rm div}u-P(\rho),\\
&L=\frac{1}{\rho}(\p_{y}A+\p_{x}B)\;\;\;\mbox{and}\;\;\;H=\frac{1}{\rho}(-\p_{x}A+\p_{y}B).
\end{aligned}
$$
%In fact we have :
%$$
%\begin{aligned}
%2{\rm div}(\mu Du)+\n(\lambda(\rho){\rm div}u)&=\mu\D u+\n\big((\mu+\lambda(\rho)){\rm div}u\big)\\
%&=\mu{\rm div}({\rm curl}u)+\n\big((2\mu+\lambda(\rho)){\rm div}u\big)
%\end{aligned}
%$$
%\\
%\\
%\\
We now want to obtain some estimates on the unknowns $A$ and $B$, let us start with rewriting the momentum equation under the following eulerian form:
\begin{equation}
\p_{t}u+u\cdot\n u-\frac{1}{\rho}\D u-\frac{1}{\rho}\n((\mu+\lambda(\rho)){\rm div}u)+\n(\frac{P(\rho)}{\gamma\rho})=0
\label{euler}
\end{equation}
Next if we apply the operator ${\rm curl}$, we get:
\begin{equation}
\p_{t}A+u\cdot\n A+A{\rm div}u=\p_{y}L-\p_{x}H.
\label{14curl}
\end{equation}
Next we apply the operator ${\rm div}$ to the momentum equation (\ref{euler}):
\begin{equation}
\p_{t}{\rm div}u+u\cdot\n u-\frac{1}{\rho}\D u-\frac{1}{\rho}\n((\mu+\lambda(\rho)){\rm div}u)+\n(\frac{P(\rho)}{\gamma\rho})=0,
\label{euler1}
\end{equation}
and via the mass equation we have:
\begin{equation}
\begin{aligned}
&\p_{t}B+U\cdot\n B-\rho(2+\lambda)\big(B(\frac{1}{2+\lambda})^{'}
+(\frac{P}{2+\lambda})^{'}\big){\rm div}u\\
&\hspace{2cm}+(2+\lambda)(U_{x}^{2}+2U_{y}V_{x}+V_{y}^{2})=(2+\lambda)(L_{x}+H_{y}).
\end{aligned}
\label{14efficace}
\end{equation}
As in \cite{VG} multiplying the equation(\ref{14curl}) by $A$ and integrate over $\Omega$ we obtain:
\begin{equation}
\begin{aligned}
&\int_{\Omega}\frac{1}{2}\frac{d}{dt}[A^{2}]dx
+\frac{1}{2}\int_{\Omega}{\rm div}u\, A^{2}dx+\int_{\Omega}(L\p_{y}A-H\p_{x}A)dx=0.
\end{aligned}
\label{in1}
\end{equation}
Similarly multiplying the equation(\ref{14efficace}) by  $\frac{1}{2+\lambda}B$ and integrate over $\Omega$ we have:
\begin{equation}
\begin{aligned}
&\int_{\Omega}\frac{1}{2+\lambda}\frac{d}{dt}(\frac{1}{2}B^{2})dx-\frac{1}{2}\int_{\Omega}{\rm div}u\frac{B^{2}}{2+\lambda}dx-\frac{1}{2}\int_{\Omega}u\cdot\n(\frac{1}{2+\lambda})B^{2}dx\\
&-\int_{\Omega}
\rho B{\rm div}u\big(B(\frac{1}{2+\lambda})^{'}+(\frac{P}{2+\lambda})^{'}\big)dx+\int_{\Omega}B(U_{x}^{2}+2U_{y}V_{x}+V_{y}^{2}
)dx\\
&\hspace{7cm}+\int_{\Omega}(L\p_{x}B+H\p_{y}B)dx=0.
\end{aligned}
\label{in11}
\end{equation}
We recall now that:
$$\p_{t}(\frac{1}{2+\lambda})+(\frac{1}{2+\lambda})^{'}\rho{\rm div}u+\n(\frac{1}{2+\lambda})\cdot u=0.$$
By combining the previous equality and (\ref{in11}) we get:
\begin{equation}
\begin{aligned}
&\frac{1}{2}\int_{\Omega}\frac{d}{dt}(\frac{1}{2+\lambda}B^{2})dx
-\frac{1}{2}\int_{\Omega}{\rm div}uB^{2}\big(\frac{1}{2+\lambda}
-\rho(\frac{1}{2+\lambda})^{'}\big)dx\\
&-\int_{\Omega}
\rho B{\rm div}u\big(B(\frac{1}{2+\lambda})^{'}+(\frac{P}{2+\lambda})^{'}\big)dx+\int_{\Omega}B({\rm div}u)^{2}
dx\\
&+2\int_{\Omega}B(\p_{y}U\p_{x}V-\p_{x}U\p_{y}V)dx+\int_{\Omega}(L\p_{x}B+H\p_{y}B)dx=0.
\end{aligned}
\label{in2}
\end{equation}
Summing (\ref{in1}) and (\ref{in2}) we have:
%$$
%\begin{aligned}
%&\int_{\Omega}\frac{1}{2}\frac{d}{dt}(A^{2}+\frac{B^{2}}{2+\lambda})dx+
%\int_{\Omega}\frac{({\rm curl}u_{y}+B_{x})^{2}+(-{\rm curl}u_{x}+B_{y})^{2}}{\rho}\\
%&-\frac{1}{2}\int_{\Omega}B^{2}\big((\frac{1}{2+\lambda})_{t}+(\frac{U}{2+\lambda})_{x}
%+\frac{V}{2+\lambda})_{y}\big)dx+\int_{\Omega}B(U_{x}^{2}+2U_{y}V_{x}+V_{y}^{2})dx\\
%&\hspace{3cm}+\int_{\Omega}\frac{1}{2}{\rm div}u({\rm curl}u)^{2}dx-\int_{\Omega}
%\rho B{\rm div}u\big(B(\frac{1}{2+\lambda})^{'}+(\frac{P}{2+\lambda})^{'}\big)dx=0
%\end{aligned}
%$$
%In view of the continuity equation, we then obtain:
$$
\begin{aligned}
&\int_{\Omega}\frac{1}{2}\frac{d}{dt}[(A^{2}+\frac{B^{2}}{2+\lambda})]dx+\int_{\Omega}\frac{1}{2}{\rm div}uA^{2}dx
+\int_{\Omega}\frac{(A_{y}+B_{x})^{2}+(-A_{x}+B_{y})^{2}}{\rho}dx\\
&-\frac{1}{2}\int_{\Omega}B^{2}{\rm div}u\big(\frac{1}{2+\lambda}-\rho(\frac{1}{2+\lambda})^{'}\big)dx
+2\int_{\Omega}B(U_{y}V_{x}-U_{x}V_{y})dx\\
&-\int_{\Omega}
\rho B{\rm div}u\big(B(\frac{1}{2+\lambda})^{'}+(\frac{P}{2+\lambda})^{'}\big)dx+
\int_{\Omega}B{\rm div}u^{2}dx=0.
\end{aligned}
$$
As:
$${\rm div}u^{2}={\rm div}u(\frac{B}{2+\lambda}+\frac{P}{2+\lambda}),$$
we deduce:
\begin{equation}
\begin{aligned}
&\int_{\Omega}\frac{1}{2}\frac{d}{dt}[(A^{2}+\frac{B^{2}}{2+\lambda})]dx+\int_{\Omega}\frac{1}{2}{\rm div}u A^{2}dx+
\int_{\Omega}\frac{(A_{y}+B_{x})^{2}+(-A_{x}+B_{y})^{2}}{\rho}dx\\
&+\int_{\Omega}\frac{1}{2}B^{2}{\rm div}u\big(\frac{1}{2+\lambda}-\rho(\frac{1}{2+\lambda})^{'}\big)dx+
\int_{\Omega}B{\rm div}u\big(\frac{P}{2+\lambda})-\rho(\frac{P}{2+\lambda}))^{'}\big)dx
\\
&+2\int_{\Omega}B(U_{y}V_{x}-U_{x}V_{y})dx-\int_{\Omega}
\rho B{\rm div}u\big(B(\frac{1}{2+\lambda})^{'}+(\frac{P}{2+\lambda})^{'}\big)dx=0.
\end{aligned}
\label{37}
\end{equation}
Let us set:
\begin{equation}
\begin{aligned}
&Z(t)=\big(\int_{\Omega}(A^{2}+\frac{B^{2}}{2+\lambda})dx\big)^{\frac{1}{2}}\\
&a(t)=\big(\int_{\Omega}\frac{(A_{y}+B_{x})^{2}+(-A_{x}+B_{y})^{2}}{\rho}dx\big)^{\frac{1}{2}}
,\;\;\;t\in[0,T].
\end{aligned}
\label{38}
\end{equation}
Next we have:
$$\int_{\Omega}\big((A_{y}+B_{x})^{2}+(-A_{x}+B_{y})^{2}\big)dx=\int_{\Omega}(
A_{x}^{2}+A_{y}^{2}+B_{x}^{2}+B_{y}^{2})dx.$$
Let us observe that for every $r$, $1\geq 4r>0$, from the result on elliptic system and  by H\"older inequalities we get as in \cite{VG}:
$$
\begin{aligned}
&\|\n A\|_{L^{2(1-r)}(\Omega)}+\|\n B\|_{L^{2(1-r)}(\Omega)}\leq C
\big(\int_{\Omega}\frac{(A_{y}+B_{x})^{2}+(-A_{x}+B_{y})^{2}}{\rho}dx\big)^{\frac{1}{2}}\\
&\hspace{10cm}\big(\int_{\Omega}\rho^{\frac{1-r}{r}}dx\big)^{\frac{r}{2(1-r)}}.
\end{aligned}
$$
From (\ref{36}), we have:
\begin{equation}
\big(\|\n A\|_{L^{2(1-r)}(\Omega)}+\|\n B\|_{L^{2(1-r)}(\Omega)}\big)\leq C(\frac{1}{r})^{\frac{1}{2(\beta-1)}}a(t)
\label{39}
\end{equation}
\begin{remarka}
Let us point out that the estimate (\ref{39}) is better than the corresponding one in \cite{VG} due to the better estimate (\ref{36}).
\end{remarka}
Moreover via (\ref{36}) we also obtain the following inequality:
\begin{equation}
\big(\|\n u\|_{L^{2}(\Omega)}+\|A\|_{L^{2}(\Omega)}+\|\sqrt{2+\lambda}\,{\rm div}u\|_{L^{2}(\Omega)}\big)\leq
C(1+Z(t)),\;\;\;t\in[0,T].
\label{40}
\end{equation}
Now, are interested in providing other estimates for the non positive terms of the equality (\ref{37}).
\subsubsection*{Estimates for the terms of (\ref{37})}
Following \cite{VG}, using (\ref{17}), the lemma \ref{lemme1} (with $\alpha=\frac{1-\e}{2(1-2\e)}$) and Young's inequality (with $p=\frac{2(1-2\e)}{1-\e}$, $q=\frac{2(1-2\e)}{1-3\e}$ and $p_{1}=\frac{(1-2\e)(2+\e)}{1-\e}$, $q_{1}=\frac{(1-2\e)(2+\e)}{1-2\e-2\e^{2}}$) the , we obtain:
$$
\begin{aligned}
\big|\frac{1}{2}\int_{\Omega}{\rm div}u\,A^{2}\big|\leq& \frac{1}{2}\|{\rm div}u(t)\|
_{L^{2}}\|A(t)\|_{L^{4}}^{2},\\
&\leq C \|{\rm div}u\|_{L^{2}}(^{\frac{1}{2}}\|A\|_{L^{2}})^{\frac{1-3\e}{1-2\e}}(^{\frac{1}{2}}\|\n A\|_{L^{2(1-\e)}})^{\frac{1-\e}{1-2\e}}\\
&\leq C(1+Z(t))Z(t)^{\frac{1-3\e}{1-2\e}}\big((\frac{1}{\e})^{\frac{1}{2(\beta-1)}}a(t)\big)^{\frac{1-\e}{1-2\e}}\\
&\leq \delta a^{2}(t)+C(\delta)(1+Z(t))^{\frac{2(1-2\e)}{1-3\e}}Z(t)^{2}(\frac{1}{\e})^{\frac{1}{\beta-1}}\\
&\leq \delta a^{2}(t)+C(\delta)(1+Z(t)^{2})^{2+\frac{\e}{1-3\e}}(\frac{1}{\e})^{\frac{1}{\beta-1}}.
\end{aligned}
$$
%Une idee ici peut-etre de traiter le cas des petits temps via des solutions fortes et une donnee vitesse legerement surcritique ce qui permettrait un ${\rm div}u\in L^{2}(L^{p})$ et ainsi pouvoir utiliser un Gronwall classique).\\
We  now are interested in estimating the term in (\ref{37}) corresponding to:
$$
\begin{aligned}
&I_{1}=\big|\frac{1}{2}\int_{\Omega}B^{2}{\rm div}u\big(\frac{1}{2+\lambda}-\rho(\frac{1}{2+\lambda})^{'}\big)dx\big|\\
&\hspace{4,5cm}=\big|\frac{1}{2}\int_{\Omega}B^{2}(\frac{B}{2+\lambda}+\frac{P}{2+\lambda})\big(\frac{1}{2+\lambda}-\rho(\frac{1}
{2+\lambda})^{'}\big)dx\big|
\end{aligned}
$$
Easily  there exist a positive constant $C>0$ such that:
$$
\begin{aligned}
\big|\frac{1}{2+\lambda}-\rho(\frac{1}{2+\lambda})^{'}\big|&\leq C,
\end{aligned}
$$
 for all $\rho\in[0,+\infty)$. We deduce that:
$$I_{1}\leq C (m^{'})^{\beta}\big(\int_{\Omega}\frac{|B|^{3}}{2+\lambda}dx+\int_{\Omega}\frac{|B|^{2}}{2+\lambda}|P|dx\big).$$
By Young's inequality we have:
$$
\begin{aligned}
\big|\int_{\Omega}B{\rm div}u\big(\frac{P}{2+\lambda}-\rho&(\frac{P}{2+\lambda})^{'}\big)dx\big|=
\big|\int_{\Omega}B(\frac{B}{2+\lambda}+\frac{P}{2+\lambda})\big(\frac{P}{2+\lambda}-\rho(\frac{P}{2
+\lambda})^{'}\big)dx\big|\\
&\hspace{2cm}\leq C (1+\int_{\Omega}\frac{|B|^{3}}{2+\lambda}dx).
\end{aligned}
$$
Now, the last term in (\ref{37}) can be treated as follows:
$$\big|2\int_{\Omega}B(U_{y}V_{x}-U_{x}V_{y})dx\big|\leq\int_{\Omega}|B|(U_{x}^{2}+U_{y}^{2}+V_{x}^{2}+V_{y}^{2})dx$$
Via the previous estimate, the notations (\ref{38}), and using the equality (\ref{37}),
we get:
\begin{equation}
\begin{aligned}
&\frac{1}{2}(Z^{2}(t))+a^{2}(t)\leq\delta a^{2}(t)+C(\delta)(1+Z(t)^{2})^{1+\frac{\e}{1-3\e}}(\frac{1}{\e})^{\frac{2}{\beta-1}}\\
&+C(1+\int_{\Omega}\frac{|B|^{3}}{2+\lambda}dx)+\int_{\Omega}|B|(U_{x}^{2}+U_{y}^{2}+V_{x}^{2}+V_{y}^{2})dx.
\end{aligned}
\label{41}
\end{equation}
It remains to estimate the two last terms on the right hand side of (\ref{41}). In this goal, from (\ref{9}) we have:
\begin{equation}
\|B\|_{L^{2m}(\Omega)}\leq C(\|B\|_{L^{1}(\Omega)}+m^{\frac{1}{2}}\|\n B\|^{1-s}_{L^{\frac{2m}{m+\e}}(\Omega)}
\|B\|^{s}_{L^{2(1-\e)}(\Omega)})
\label{42}
\end{equation}
where: $s=\frac{(1-\e)^{2}}{m-\e(1-\e)}$ and $C>0$ is a positive constant independent of $m>2$.\\
Now in inequalities (\ref{41}) and (\ref{42}) we fix $\e=2^{-m}$ with $m>2$. Using estimate (\ref{36}) for the density, we derive the inequalities:
$$
\begin{aligned}
\|B\|_{L^{1}(\Omega)}=\int_{\Omega}|B|dx&=\int_{\Omega}(\frac{1}{2+\lambda})^{\frac{1}{2}}|B|(2+\lambda)^{\frac{1}{2}}dx
\leq \|(2+\lambda)^{\frac{1}{2}}\|_{L^{2}(\Omega)}Z(t)\\
&\leq C Z(t).
\end{aligned}
$$
$$
\begin{aligned}
&\|B\|^{s}_{L^{2(1-\e)}(\Omega)}=\big(\int_{\Omega}(\frac{1}{2+\lambda})
^{1-\e}|B|^{2(1-\e)}(2+\lambda)^{1-\e}dx\big)^{\frac{s}{2(1-\e)}}\\
&\hspace{1cm}\leq Z(t)^{s}\|2+\lambda\|_{L^{\frac{1-\e}{\e}}(\Omega)}^{\frac{s}{2}}\leq C(\frac{1}{\e})^{\frac{\beta s}{2(\beta-1)}}Z^{s}(t)\leq C 2^{2sm}Z^{s}(t
),\\
&\hspace{1cm}\leq C Z^{s}(t).
\end{aligned}
$$
From inequalities (\ref{39}) and (\ref{42}) we finally obtain:
\begin{equation}
\begin{aligned}
\|B\|_{L^{2m}(\Omega)}&\leq C\big(Z(t)+m^{\frac{1}{2}}(\frac{m}{\e})^{\frac{1-s}{2(\beta-1)}}(a(t))^{1-s} Z^{s}(t)\big),\\
&\leq C\big(Z(t)+m^{\frac{1}{2}}(\frac{m}{\e})^{\frac{1-s}{2(\beta-1)}}(a(t))^{1-s} Z^{s}(t)\big).
\end{aligned}
\label{43}
\end{equation}
Now dealing with the integral with $|B|^{3}$, we have:
$$
\begin{aligned}
&\int_{\Omega}\frac{|B|^{3}}{2+\lambda}dx=\int_{\Omega}\frac{|B|^{2-\frac{1}{m-1}}}{(2+\lambda)^{1-\frac{1}{2(m-1)}}}
(\frac{1}{2+\lambda})^{\frac{1}{2(m-1)}}|B|^{1+\frac{1}{m-1}}dx\\
&\leq\int_{\Omega}\frac{|B|^{2-\frac{1}{m-1}}}{(2+\lambda)^{1-\frac{1}{2(m-1)}}}
|B|^{\frac{m}{m-1}}dx\leq Z(t)^{2-\frac{1}{m-1}}\|B\|^{\frac{m}{m-1}}_{L^{2m}(\Omega)},\\
&\leq(\int_{\Omega}\frac{|B|^{2}}{2+\lambda}dx)^{1-\frac{1}{2(m-1)}}(\int_{\Omega}
|B|^{2m}dx)^{\frac{1}{2(m-1)}}\leq ^{-1+\frac{1}{2(m-1)}}Z(t)^{2-\frac{1}{m-1}}\|B\|^{\frac{m}{m-1}}_{L^{2m}(\Omega)}\\
&\leq C Z^{2-\frac{1}{m-1}}(t)\big(Z(t)^{\frac{m}{m-1}}+m^{\frac{m}{2(m-1)}}
(\frac{m}{\e})^{\frac{m(1-s)}{2(\beta-1)(m-1)}}(a(t))^{\frac{m(1-s)}{m-1}} Z^{\frac{ms}{m-1}}(t)\big),\\
&\leq C \big(Z(t)^{3}+m^{\frac{m}{2(m-1)}}
(\frac{m}{\e})^{\frac{m(1-s)}{2(\beta-1)(m-1)}}(a(t))^{\frac{m(1-s)}{m-1}} Z^{2+\frac{ms-1}{m-1}}(t)\big),
\end{aligned}
$$
where $C>0$ is a positive constant independent of $m>2$ and $\e=2^{-m}$. Finally applying applying Young's inequality with $p=\frac{2(m-1)}{m(1-s)}$ and $q=\frac{2(m-1)}{m(s+1)-2}$ we have:
\begin{equation}
\begin{aligned}
&\int_{\Omega}\frac{|B|^{3}}{2+\lambda}dx\leq
C \big(Z^{3}(t)+m^{\frac{1}{2}}(\frac{m}{\e})^{\frac{1}{2(\beta-1)}}
a^{\frac{m(1-s)}{m-1}}(t)Z^{2+\frac{ms-1}{m-1}}(t)\big),\\
&\leq\delta a^{2}(t)+C(\delta)\big(Z^{3}(t)+m^{\frac{m-1}{m(s+1)-2}}(\frac{m}{\e})^{\frac{m-1}{(\beta-1)(m(s+1)-2)}}
Z^{4+\frac{2(1-ms)}{m(s+1)-2}}(t)\big)\\
\end{aligned}
\label{44}
\end{equation}
From (\ref{42}), we verify that:
\begin{equation}
1-ms=1-\frac{m(1-\e)^{2}}{m-\e(1-\e)}=\e\frac{m(2-\e)+\e-1}{m-\e(1-\e)},
\label{45}
\end{equation}
and then:
$$\lim_{m\rightarrow+\infty}(2^{m}(1-ms))=2,$$
hence via (\ref{44}) and (\ref{45}) we get:
\begin{equation}
\int_{\Omega}\frac{|B|^{3}}{2+\lambda}dx\leq\delta a^{2}(t)+C(\delta)\big((1+Z^{2}(t))^{2}+m(\frac{m}{\e})^{\frac{1}{\beta-1}}
(1+Z^{2}(t))^{2+\frac{1-ms}{m(s+1)-2}}\big).
\label{46}
\end{equation}
Now, consider the last term in (\ref{41}):
$$I_{2}=\int_{\Omega}|B|(U_{x}^{2}+U_{y}^{2}+V_{x}^{}+V_{y}^{2})\leq\|B\|_{L^{2m}(\Omega)}(\int_{\Omega}(|\n U|^{2}+|\n V|^{2})^{\frac{2m}{2m-1}}dx)^{1-\frac{1}{2m}}.$$
Recalling the relation $3>\frac{4m}{2m-1}>2$, $m>2$, from the properties of elliptic system \cite{29} we derive the inequality:
$$(\int_{\Omega}(|\n U|^{2}+|\n V|^{2})^{\frac{2m}{2m-1}}dx)^{1-\frac{1}{2m}}\leq C(\|{\rm div}u\|^{2}_{
L^{\frac{4m}{2m-1}}(\Omega)}+\|A\|^{2}_{
L^{\frac{4m}{2m-1}}(\Omega)}).$$
Thus the previous inequality furnish the estimate:
\begin{equation}
I_{2}\leq C\|B\|_{L^{2m}(\Omega)}(\|{\rm div}u\|^{2}_{
L^{\frac{4m}{2m-1}}(\Omega)}+\|A\|^{2}_{
L^{\frac{4m}{2m-1}}(\Omega)}).
 \label{47}
\end{equation}
Next we have as $A$ vanishes on the boundary of the domain $\Omega$, we applythe Gagliardo-Niremberg inequality:
\begin{equation}
\begin{aligned}
&\|A\|^{2}_{
L^{\frac{4m}{2m-1}}(\Omega)}\leq C\|A\|^{2-\frac{1-\e}{m(1-2\e)}}_{L^{2}(\Omega)}
\|\n A\|^{\frac{1-\e}{m(1-2\e)}}_{L^{2(1-\e)}(\Omega)},\\
&\leq C Z^{2-\frac{1-\e}{m(1-2\e)}}(t)\big((\frac{1}{\e})^{\frac{1}{2(\beta-1)}}a(t)\big)^{\frac{1-\e}{m(1-2\e)}}\leq 
C Z^{2-\frac{1-\e}{m(1-2\e)}}(t)(a(t))^{\frac{1-\e}{m(1-2\e)}}.
\end{aligned}
 \label{48}
\end{equation}
Since $B=(2+\lambda){\rm div}u-P$, estimate (\ref{36}) provides:
\begin{equation}
 \begin{aligned}
&\|{\rm div}u\|^{2}_{
L^{\frac{4m}{2m-1}}(\Omega)}=\|\frac{B}{2+\lambda}+\frac{P}{2+\lambda}\|^{2}_{
L^{\frac{4m}{2m-1}}(\Omega)},\\
&\leq C( \|\frac{B}{2+\lambda}\|^{2}_{
L^{\frac{4m}{2m-1}}(\Omega)}+1).
 \end{aligned}
\label{49}
\end{equation}
We can now deal with the right-hand side of (\ref{49}) as follows:
$$
\begin{aligned}
& \|\frac{B}{2+\lambda}+\frac{P}{2+\lambda}\|^{2}_{
L^{\frac{4m}{2m-1}}(\Omega)}\leq\big(\int_{\Omega}\frac{|B|^{\frac{2m(2m-3)}{(m-1)(2m-1)}}}{2+\lambda}|B|^{\frac{2m}{(m-1)(2m-1)}}dx\big)^{1-\frac{1}{2m}},\\
&
\leq\|B\|^{\frac{1}{m-1}}_{L^{2m}(\Omega)}\big(\int_{\Omega}\frac{|B|^{2}}{(2+\lambda)^{\frac{(m-1)(2m-1)}{m(2m-3)}}}dx
\big)^{\frac{2m-3}{2m-2}}\leq\|B\|^{\frac{1}{m-1}}_{L^{2m}(\Omega)}\big(\int_{\Omega}\frac{|B|^{2}}{(2+\lambda)}dx
\big)^{\frac{2m-3}{2m-2}},\\
&\leq\|B\|^{\frac{1}{m-1}}_{L^{2m}(\Omega)}(Z(t))^{\frac{2m-3}{2m-2}}=\|B\|^{\frac{1}{m-1}}_{L^{2m}(\Omega)}
(Z(t))^{2-\frac{1}{m-1}}.
\end{aligned}
$$
Thus,
\begin{equation}
\|{\rm div}u\|^{2}_{
L^{\frac{4m}{2m-1}}(\Omega)}\leq C(1+(Z(t))^{2-\frac{1}{m-1}}\|B\|^{\frac{1}{m-1}}_{L^{2m}(\Omega)}).
 \label{50}
\end{equation}
Using estimates (\ref{48}) and (\ref{50}), from (\ref{47}) we have:
$$I_{2}\leq C\|B\|_{L^{2m}(\Omega)}(1+Z^{2-\frac{1-\e}{m(1-2\e)}}(t)(a(t))^{\frac{1-\e}{m(1-2\e)}}+
Z^{2-\frac{1}{m-1}}(t)\|B\|_{L^{2m}(\Omega)}^{\frac{1}{m-1}}\big).$$
Using estimate (\ref{43}) for $\|B\|_{L^{2m}(\Omega)}$, we finally get:
\begin{equation}
 \begin{aligned}
&I_{2}\leq C\big((Z(t))^{3-\frac{1-\e}{m(1-2\e)}}(a(t))^{\frac{1-\e}{m(1-2\e)}}+Z(t)+Z^{3}(t)\\
&+m^{\frac{1}{2}}
(\frac{m}{\e})^{\frac{1-s}{2(\beta-1)}}Z^{2+s-\frac{1-\e}{m(1-2\e)}}(t)(a(t))^{1-s+\frac{1-\e}{m(1-2\e)}}\\
&+m^{\frac{1}{2}}
(\frac{m}{\e})^{\frac{m(1-s)}{2(\beta-1)(m-1)}}Z^{2+\frac{ms-1}{m-1}}(t)(a(t))^{\frac{m(1-s)}{m-1}}
+m^{\frac{1}{2}}
(\frac{m}{\e})^{\frac{1-s}{2(\beta-1)}}Z^{s}(t)(a(t))^{1-s}\big).
 \end{aligned}
\label{51}
\end{equation}
Using Young's inequality, we treat the summand in (\ref{51}) as follows:
$$
\begin{aligned}
& C(Z(t))^{3-\frac{1-\e}{m(1-2\e)}}(a(t))^{\frac{1-\e}{m(1-2\e)}}\leq\delta a^{2}(t)+C(1+Z^{2}(t))^{2},\\
&C(Z(t)+Z^{3}(t))\leq C(1+Z^{2}(t))^{2},\\[3mm]
&Cm^{\frac{1}{2}}
(\frac{m}{\e})^{\frac{1-s}{2(\beta-1)}}Z^{2+s-\frac{1-\e}{m(1-2\e)}}(t)(a(t))^{1-s+\frac{1-\e}{m(1-2\e)}}\leq\\
&\hspace{1cm}\delta a^{2}(t)+Cm(\frac{m}{\e})^{\frac{1}{\beta-1}}(1+Z^{2}(t))^{2+\frac{1-ms+(2ms-1)\e}{(1+s)m(1-2\e)-1+\e}},\\[3mm]
&Cm^{\frac{1}{2}}
(\frac{m}{\e})^{\frac{1-s}{2(\beta-1)}}Z^{s}(t)(a(t))^{1-s}\leq\delta a^{2}(t)+Cm(\frac{m}{\e})^{\frac{1}{\beta-1}}(1+Z^{2}(t)),\\[3mm]
&m^{\frac{1}{2}}
(\frac{m}{\e})^{\frac{m(1-s)}{2(\beta-1)(m-1)}}Z^{2+\frac{ms-1}{m-1}}(t)(a(t))^{\frac{m(1-s)}{m-1}}\leq\\
&\hspace{1cm}\delta a^{2}(t)+Cm(\frac{m}{\e})^{\frac{1}{\beta-1}}(1+Z^{2}(t))^{2+\frac{1-ms}{(1+s)m-2}} 
\end{aligned}
$$
Here $\delta$ is a small positive constant to be mentioned below. From inequality (\ref{51}) we derive that:
\begin{equation}
\begin{aligned}
&I_{2}\leq \delta a^{2}(t)+C\big(m(\frac{m}{\e})^{\frac{1}{\beta-1}}(1+Z^{2}(t))^{2+\frac{1-ms+(2ms-1)\e}{(1+s)m(1-2\e)-1+\e}}\\
&+(1+Z^{2}(t))^{2}+m(\frac{m}{\e})^{\frac{1}{\beta-1}}(1+Z^{2}(t))+m(\frac{m}{\e})^{\frac{1}{\beta-1}}(1+Z^{2}(t))^{2+\frac{1-ms}{(1+s)m-2}} 
\end{aligned}
 \label{52}
\end{equation}
From (\ref{46}) and (\ref{52}), and inequality (\ref{41}) we have:
\begin{equation}
\begin{aligned}
&\frac{1}{2}\frac{d}{dt}(Z^{2}(t)+a^{2}(t))\leq\delta a^{2}(t)+C(\delta)(1+Z^{2}(t))^{2+\frac{\e}{3-\e}}(\frac{1}{\e})^{\frac{1}{\beta-1}}\\
&+C(1+\delta a^{2}(t))+C C(\delta)(1+Z^{2}(t))^{2}+C C(\delta)\big( m(\frac{m}{\e})^{\frac{1}{\beta-1}}(1+Z^{2}(t))^{2+\frac{1-ms}{m(s+1)-2}}\\
&+4\delta a^{2}(t)+C\big(((1+Z^{2}(t))^{2}+m(\frac{m}{\e})^{\frac{1}{\beta-1}}(1+Z^{2}(t))^{2+\frac{1-ms
+(2ms-1)\e}{(1+s)m(1-2\e)-1+\e}}\\
&+m(\frac{m}{\e})^{\frac{1}{\beta-1}}(1+Z^{2}(t))^{2}+m(\frac{m}{\e})^{\frac{1}{\beta-1}}
(1+Z^{2}(t))^{2+\frac{1-ms}{m(s+1)-2}}\big).
\end{aligned}
\label{53}
\end{equation}
Choose $\delta>0$ such that:
$$5\delta+\delta C=\frac{1}{2}.$$
Since $s=\frac{(1-\e)^{2}}{m-\e(1-\e)}$ and $\e=2^{-m}$, $m>2$, we have:
$$
\begin{aligned}
&\frac{1-ms}{m(s+1)-2}\leq 4\e,\;\;\frac{1-ms
+(2ms-1)\e}{(1+s)m(1-2\e)-1+\e}\leq 4\e\;\;\;\mbox{and}\;\;\;\frac{\e}{1-3\e}\leq 4\e.
\end{aligned}
$$
Then by (\ref{53}) and the fact that $Z^{2}(t)\in L^{1}(0,T)$, we obtain the inequality with $0<\bar{T}<\frac{T}{2}$:
\begin{equation}
\frac{1}{2}\frac{d}{dt}(1+Z^{2}(t))+a^{2}(t)\leq m(\frac{m}{\e})^{\frac{1}{\beta-1}}
(1+Z^{2}(t))^{2+4\e},
\label{54}
\end{equation}
From (\ref{54}) we have for $0\leq t<T$:
$$
\frac{1}{(1+Z^{2}(t))^{4\e}}-\frac{1}{(1+Z^{2}(\bar{T}))^{4\e}}+Cm\e(\frac{m}{\e})^{\frac{1}{\beta-1}}
\geq0.$$
\begin{remarka}
Let us point out that the last inequality is better than in \cite{VG} and allows us to assume only $\beta>2$.
\end{remarka}
Now, take $N>2$ such that:
$$1-CN\e(\frac{N}{\e})^{\frac{1}{\beta-1}}(1+Z^{2}(\bar{0}))^{4\e}\geq\frac{1}{2},\;\e=2^{-N}.$$
Here the fact that $\beta>2$ allows to conclude and by this fact improve the results of \cite{VG}.
We get finally that for $0\leq t<T$:
\begin{equation}
Z^{2}(t)\leq 2^{2^{N-2}}(1+Z^{2}(0))-1,\;\;\;t\in[0,T].
\label{55}
\end{equation}
Now, from inequality (\ref{54}) we get moreover that:
\begin{equation}
\int^{T}_{0}a^{2}(t)dt\leq C.
\label{56}
\end{equation}
Now by estimate (\ref{36}) for the density, there exists a positive constant $C$ depending continuously on the data of the problem and such that:
\begin{equation}
\begin{aligned}
&\sup_{0<t<T}\int_{\Omega}\biggl(({\rm curl}u)^{2}+\frac{1}{2+\lambda(\rho)}\big((2+\lambda(\rho)){\rm div}u-P(\rho)\big)^{2}\biggl)(t,x)dx\leq C,\\
&\sup_{0<t<T}\int_{\Omega}\big(({\rm curl}u)^{2}+(2+\lambda(\rho))({\rm div}u)^{2}\big)(t,x)dx\leq C,\\
&\int^{T}_{0}\int_{\Omega}\frac{(A_{y}+B_{x})^{2}+(-A_{x}+B_{y})^{2}}{\rho}dxdy\leq C.
\end{aligned}
\label{57}
\end{equation}
The rest of the proof follows exactly the same lines than in \cite{VG} and then we refer to \cite{VG}.
%\\
%Il est a signaler que si l'on choisit $\beta=1+\frac{2}{m}$ alors avec m tres grand on a:
%$$\|\rho\|_{L^{2m\beta+1}}\leq C m^{\frac{-3\beta+2}{2\beta-1}},\;\;\;t\in[0,T],$$
%qui est par consequent tres petit. Que se passe-t-il si on construit une suite de $\rho_{n}$, alors la limite tend a l'infini vers une fonction $\rho$ $L^{\infty}$.
\section{Proof of theorem \ref{theo2}}
\label{section6}
The existence part of the theorem is proved by an iterative method. We define a sequence $(q^{n},m^{n})$ such that:
$$
\begin{cases}
\begin{aligned}
&\p_{t}q^{0}-\mu\D q^{0}+{\rm div}m^{0}=0,\\
&\p_{t}m^{0}-\mu\D m^{0}+r m^{0}=0,\\
&(q^{0},m^{0})=(q_{0},m_{0}).
\end{aligned}
\end{cases}
$$
Assuming that $(q^{n},m^{n})$ is in $E_{T}$ with:
$$E_{T}=\big(\widetilde{C}_{T}(B^{\N}_{2,1})\cap L^{1}_{T}(B^{\N+2}_{2,1})\big)\times\big(
\widetilde{C}_{T}(B^{\N-1}_{2,1})\cap L^{1}_{T}(B^{\N+1}_{2,1}\cap B^{\N-1}_{2,1})\big)^{N},$$
we define then  $q^{n+1}=q^{0}+\bar{q}^{n+1}$, $m^{n+1}=m^{0}+\bar{m}^{n+1}$  such that 
$(\bar{q}_{n+1},\bar{m}_{n+1})$ be the solution
of the following system:
$$
\begin{cases}
\begin{aligned}
&\p_{t}\bar{q}^{n+1}-\mu\D\bar{q}^{n+1}+{\rm div}\bar{m}^{n+1}=0,\\
& \p_{t}\bar{m}^{n+1}-\mu\D\bar{m}^{n+1}+r \bar{m}^{n+1} =G^{n},\\
&(\bar{q}^{n+1},\bar{m}^{n+1})_{/t=0}=(0,0),
\end{aligned}
\end{cases}
$$
with:
$$
\begin{aligned}
%F^{n}=&-{\rm div}(m^{n})\\
G^{n}=&-{\rm div}(\frac{m^{n}}{h^{n}}\otimes m^{n})
\end{aligned}
$$
We also set: $h^{n}=q^{n}+1$.
\subsubsection*{1) First step, uniform bounds:}
 Let $\e$ be a small
parameter and by proposition \ref{5chaleur}, we have for any $T>0$:
\begin{equation}
\begin{aligned}
&\|q^{0}\|_{L^{\infty}_{T}(B^{\N}_{2,1})\cap L^{1}_{T}(B^{\N+2}_{2,1})}\leq C\|q_{0}\|_{B^{\N}_{2,1}},\\
&\|m^{0}\|_{L^{\infty}_{T}(B^{\N-1}_{2,1})\cap L^{1}_{T}(B^{\N-1}_{2,1}\cap B^{\N+1}_{2,1})}\leq C\|m_{0}\|_{B^{\N-1}_{2,1}}.
\end{aligned}
\label{initial}
\end{equation}
We are going to show by induction that for $\e>0$ small enough:
$$\|(\bar{q}^{n},\bar{m}^{n})\|_{F_{T}}\leq\e.\leqno{({\cal{P}}_{n})}$$
As $(\bar{q}_{0},\bar{m}_{0})=(0,0)$ the result
is true for $n=0$. We suppose now $({\cal{P}}_{n})$ true and we are
going to show $({\cal{P}}_{n+1})$.
\\
To begin with we are going to show that $1+q^{n}$ is  positive. Indeed we have: $h^{0}=h^{0}_{1}+h^{0}_{2}$ such that:
$$
\begin{aligned}
&\p_{t}h^{0}_{1}-\mu \D h^{0}_{1}=0,\\
&(h^{0}_{1})_{/t=0}=h_{0}.
\end{aligned}
$$
and:
$$
\begin{aligned}
&\p_{t}h^{0}_{2}-\mu \D h^{0}_{2}=-{\rm div}m^{0},\\
&(h^{0}_{1})_{/t=0}=h_{0}.
\end{aligned}
$$
By proposition (\ref{5chaleur}) and (\ref{initial}) we have for any $T>0$:
\begin{equation}
\|h^{0}_{2}\|_{L^{\infty}_{T}(B^{\N}_{2,1})}\leq C\|m_{0}\|_{B^{\N-1}_{2,1}}.
\label{ss5}
\end{equation}
By maximum principle, we have for any $t>0$:
$$h^{0}_{1}(t,x)\geq \min_{x\in\R^{N}}h_{0}(x)\geq c>0.$$
We deduce that for $\eta=\|m_{0}\|_{B^{\N-1}_{2,1}}$ (at least inferior to $\frac{c}{4C}$ with the $C$ of (\ref{ss5})) small enough and any $t>0$:
$$h^{0}(t,x)/geq \frac{3c}{4}>0,$$
and
$$q^{0}(t,x)\geq \frac{3c}{4}-1.$$
and by definition of $q^{n}$ and the assumption ${\cap P}_{n}$ that:
$$q^{n}(t,x)\geq  \frac{3c}{4}-1-\e.$$
In particular for $\e$ small enough at least $\e\leq\frac{c}{4}$, we deduce that:
\begin{equation}
h^{n}=1+q^{n}\geq \frac{c}{2}>0.
\label{vide}
\end{equation}
In order to bound $(\bar{q}^{n},\bar{m}^{n})$ in $E_{T}$, we shall use proposition \ref{5chaleur} and in particular estimating %$F^{n}$ and 
$G^{n}$ %respectively in %$L^{1}_{T}(B^{\N}_{2,1})$ and 
in $L^{1}_{T}(B^{\N-1}_{2,1})$.
%According to proposition
%\ref{produit1},
%we have:
%\begin{equation}
%\begin{aligned}
%&\|{\rm div}(m^{n})\|_{
%L^{1}_{T}(B^{\N}_{2,1})}\leq \|m^{0}\|_{
%L^{1}_{T}(B^{\N+1}_{2,1})}+\|\bar{m}^{n}\|_{
%L^{1}_{T}(B^{\N+1}_{2,1})}.
%\end{aligned}
%\label{a1}
%\end{equation}
%Let us estimate $G^{n}$ in $L_{T}^{1}(B^{\N-1}_{2,1})$. 
By using proposition \ref{produit1}, (\ref{vide}) and lemma\ref{composition}, we obtain:
\begin{equation}
\begin{aligned}
&\|{\rm div}(\frac{m^{n}}{h^{n}}\otimes m^{n})\|_{L^{1}_{T}(B^{\N-1}_{2,1})}\leq \|\frac{m^{n}}{h^{n}}\otimes m^{n}\|_{L^{1}_{T}(B^{\N-1}_{2,1})},\\
&\leq C\|m^{n}\|_{L^{2}_{T}(B^{\N}_{2,1})}^{2}\big(\|\frac{1}{1+q^{n}}-1\|_{L^{\infty}_{T}(B^{\N}_{2,1})}+1\big),\\
&\leq C(\|m^{0}\|_{L^{2}_{T}(B^{\N}_{2,1})}^{2}+\|\bar{m}^{n}\|_{L^{2}_{T}(B^{\N}_{2,1})})^{2}(1+C(\|\frac{1}{h^{n}}\|_{L_{T}^{\infty}}))\big(\|q^{n}\|_{L_{T}^{\infty}(B^{\N}_{2,1})}+1\big),\\
&\leq C(\|m^{0}\|_{L^{2}_{T}(B^{\N}_{2,1})}^{2}+\|\bar{m}^{n}\|_{L^{2}_{T}(B^{\N}_{2,1})})^{2}\big(\|\bar{q}^{n}\|_{L_{T}^{\infty}(B^{\N}_{2,1})}+\|q^{0}\|_{L_{T}^{\infty}(B^{\N}_{2,1})}
+1\big),
\end{aligned}
\label{a2}
\end{equation}
Therefore by using %(\ref{a1}), 
(\ref{a2}), the proposition \ref{5chaleur} and $({\cal P}_{n})$ we obtain for any $T>0$:
\begin{equation}
\begin{aligned}
\|(\bar{q}^{n+1},\bar{m}^{n+1})\|_{F_{T}}&\leq C (\|m^{0}\|_{L^{2}_{T}(B^{\N}_{2,1})}^{2}+\e)^{2}\big(\e+\|q^{0}\|_{L_{T}^{\infty}(B^{\N}_{2,1})}
+1\big),\\
&\leq C (\eta+\e)^{2}\big(2+\|q^{0}\|_{L_{T}^{\infty}(B^{\N}_{2,1})}\big)
\end{aligned}
\end{equation}
By choosing $\eta=\e$ and $\e\leq \frac{1}{2C(2+\|q^{0}\|_{L_{T}^{\infty}(B^{\N}_{2,1})})}$, this implies $({\cal P})_{n+1}$.  We
have shown by induction that $(q^{n},m^{n})$ is uniformly bounded
in $F_{T}$ for any $T>0$.
\subsubsection*{Second Step: Convergence of the
sequence}
 We shall prove
that $(q^{n},m^{n})$ is a Cauchy sequence in the Banach
space $F_{T}$, hence converges to some
$(q,m)\in F_{T}$.\\
Let:
$$\delta q^{n}=q^{n+1}-q^{n}\;\;\mbox{and}\;\;\delta m^{n}=m^{n+1}-m^{n}.$$
The system verified by $(\de q^{n},\de m^{n})$ reads:
$$
\begin{cases}
\begin{aligned}
&\p_{t}\delta q^{n}-\mu\D \delta q^{n}+{\rm div}\delta m^{n}=0,\\
&\p_{t}\delta m^{n}-\mu\D\delta m^{n}+r \delta m^{n} =G^{n}-G{n-1},\\
&\delta q^{n}(0)=0\;,\;\delta u^{n}(0)=0.
\end{aligned}
\end{cases}
$$
Applying propositions \ref{5chaleur} and using $({\cal{P}}_{n})$, we get for any $T>0$:
$$
\begin{aligned}
&\|(\de q^{n},\de m^{n})\|_{F_{T}}\leq\;
C\|G^{n}-G_{n-1}\|_{L^{1}_{T}(B^{N/2-1})},\\
&\leq C\big(\|\frac{\delta m^{n}}{h^{n}}\otimes m^{n}\|_{L^{1}_{T}(B^{N/2}_{2,1})}+
\|\frac{\delta m^{n}}{h^{n}}\otimes m^{n-1}\|_{L^{1}_{T}(B^{N/2}_{2,1})}+\|m^{n}\otimes m^{n-1}(\frac{1}{h^{n}}-\frac{1}{h^{n-1}})\|_{L^{1}_{T}(B^{N/2}_{2,1})}\big).
\end{aligned}
$$
By using proposition \ref{produit1} and lemma \ref{composition}, we get:
$$\|(\de q^{n},\de m^{n})\|_{F_{T}}\leq
C\e\|(\de q^{n-1},\de m^{n-1})\|_{F_{T}}.$$ So by taking $\e$ enough small we have proved that
$(q^{n},m^{n})$ is a Cauchy sequence in $F_{T}$ which is a Banach space. It implies that $(q^{n},m^{n})$ converge to a limit $(q,m)$  in $F_{T}$.  It is easy to verify that $(q,m)$ is a
solution of the system (\ref{3system3}).
\subsubsection*{3)Uniqueness of the solution:}
The proof is similar to the proof of contraction, indeed we need the same type
 of estimates. Let us consider two solutions in
$E_{T}$: $(q_{1},m_{1})$ and
$(q_{2},m_{2})$ of the system (\ref{3system3}) with the same
initial data. With no loss of generality, one can assume that
$(q_{1},m_{1})$ is the solution found in the previous
section. We thus have:
$$q_{1}(t,x)\geq -\frac{1}{2}.
\leqno{\cal{(H)}}$$ %Let $\bar{T}$ be the largest time such that
%$q_{2}$ verifies $\cal{(H)}$. By continuity, we have $0<\bar{T}\leq
%T$.
We note:
$$\delta q=q_{2}-q_{1},\;\delta m=m_{2}-m_{1},$$ 
which verifies the system:
$$
\begin{cases}
\begin{aligned}
&\p_{t}\delta q-\mu\D\delta q+{\rm div}\delta m=0,\\
&\p_{t}\delta m-\mu\D\delta m+r \delta m=-{\rm div}(\frac{m_{1}}{h_{1}}\otimes m_{1})+
{\rm div}(\frac{m_{1}}{h_{1}}\otimes m_{1})
\end{aligned}
\end{cases}
$$
By using proposition \ref{produit1}, \ref{5chaleur} and lemma \ref{composition} on $[0,T_{1}]$
with $0<T$ we have:
$$\|(\delta q,\delta m)\|_{E_{T}}\leq A(T)\|(\delta q,\delta m)\|_{E_{T}},$$
such that for $T$ small enough $A(T)\leq\frac{1}{2}$. We thus obtain: $\de q=0$, $\de m=0$ on
$[0,T]$. And we repeat the argument in order to prove that:
$\delta q=0$, $\de m=0$ on $\R^{+}$. This conclude the proof of theorem \ref{theo1}.
\hfill {$\Box$}

\end{document}